\documentclass[11pt]{amsart}
\usepackage[]{amsmath, amsthm, amsfonts, verbatim, amssymb, tikz}
\usepackage{showkeys}

\def\cD{{\mathcal D}}

\def\cO{{\mathcal O}}

\def\cG{{\mathcal G}}

\def\bZ{{\mathbb Z}}
\def\bQ{{\mathbb Q}}
\def\bC{{\mathbb C}}

\newtheorem {theo}{Theorem}

\newtheorem {lemm}{Lemma}

\newtheorem {prop}{Proposition}

\def\wtM{\widetilde M}

\def\hM{\widehat M}

\def\bs{\bigskip}
\def\ms{\medskip}
\def\ni{\noindent}

\def\oz1{d{\overline z}^1}
\def\oz2{d{\overline z}^2}
\def\oz3{d{\overline z}^3}

\def\oI{\overline I}

\def\oz{\overline z}

\def\oGamma{\overline \Gamma}

\def\oIq1{\oI_1\cdots\oI_{q-1}}
\def\oIq2{\oI_1\cdots\oI_{q-2}}

\def\dim{{\mbox{dim}}}

\def\Aut{\mbox{Aut}}

\def\Tr{\mbox{Tr}}

\def\hs{\widehat s}
\def\hht{\widehat t}
\def\hQ{\widehat Q}

\begin{document}
\ni
\title[A surface of maximal canonical degree]
{A surface of maximal canonical degree} 
\author{ Sai-Kee Yeung}
\begin{abstract} It is known since the 70's from a paper of Beauville that the degree of the rational canonical map of a smooth projective algebraic surface of general type
is at most $36$.   Though it has been conjectured that a surface with optimal canonical degree $36$ exists, 
the highest canonical degree known earlier for a minimal surface of general type was $16$ by Persson.  The purpose of this paper is to give an affirmative answer to the conjecture
by providing an explicit surface.  %The surface is constructed from a fake projective plane first studied in Prasad and Yeung [PY].

\end{abstract}
\address[]{Mathematics Department, Purdue University, West Lafayette, IN  47907
USA} \email{yeung@math.purdue.edu}
\thanks{ 
The author was partially supported by a grant from the National Science Foundation}

\maketitle

\bigskip
\noindent

\begin{center}
{\bf 1. Introduction} 
\end{center} 

\ms
\ni{\bf 1.1.}  Let $M$ be a minimal surface of general type.  Assume that  the space of canonical sections $H^0(M,K_M)$ of $M$ is non-trivial.
Let  $N=\dim_{\bC}H^0(M,K_M)=p_g$, where $p_g$ is the geometric genus.   The canonical map is defined to
be in general a rational map $\Phi_{K_M}: M\dashrightarrow P^{N-1}$.   For a surface, this is the most natural rational mapping to study
if it is non-trivial.
  Assume that $\Phi_{K_M}$ is generically finite
with degree $d$.
It is well-known from the work of Beauville [B], that $d:=\deg\Phi_{K_M}\leqslant 36.$  We call
such degree the canonical degree of the surface, and regard it $0$ if the canonical mapping does not exist or is not
generically finite.  The following open problem is an immediate consequence of the work of [B] and is implicitly hinted there.

\bs
\ni{\bf Problem} {\it What is the optimal canonical degree of a minimal surface of general type?
Is there a minimal surface of general type with canonical degree $36$?}

\ms
Though the problem is natural and well-known, the answer remains elusive since the 70's.  The problem would be solved if
a surface of canonical degree $36$ could be constructed.    Prior to this work, the highest canonical
degree known for a surface of general type is $16$ as constructed by Persson [Pe].  We refer the readers to [DG], [Pa], [T] and [X] 
for earlier discussions on construction of surfaces with relatively large canonical degrees.  The difficulty for the open problem lies in the lack of possible candidates for
such a surface.  

Note that from the work of
[B],   a smooth surface of canonical degree $36$ is a complex two ball quotient $B_{\bC}^2/\Sigma$, where $\Sigma$ is a lattice of $PU(2,1)$.  Hence  it is infinitesimal rigid and can
neither be obtained from deformation
nor written as a complete intersection of
hypersurfaces in projective spaces. 

The purpose of this paper is to give an answer to the problem above by presenting explicitly a surface with canonical degree $36$.   Comparing to earlier methods,
we look for such a surface from a new direction, namely, arithmetic lattices coming from recent classification of fake projective planes given in [PY] and [CS].  In fact, the surface
here is constructed from  a fake projective plane 
originally studied in Prasad-Yeung [PY].

\begin{theo}
There exists a smooth minimal surface of general type $M$ with a generically finite canonical map $\Phi_{K_M}:M\rightarrow P_{\bC}^2$ of
degree $36$, constructed from an appropriate unramified covering of a well-chosen
fake projective plane of index $4$.
\end{theo}

 The example obtained above corresponds to an arithmetic lattice $\Sigma$ associated to a non-trivial division algebra over appropriate number fields as discussed in [PY].
 Arithmetic lattices coming from non-trivial division algebra are sometimes called arithmetic lattices of the second type, cf. [Ye2]. 
 In contrast, geometric complex two ball quotients studied extensively in the literature correspond to a class of examples commensurable with Deligne-Mostow surfaces in [DM], or
those constructed by Hirzebruch [H], cf.  [DM].  Further examples in this latter direction can be found in the  recent paper of Deraux-Parer-Paupert
[DPP].   The lattices involved are sometimes called arithmetic or integral lattices of first type, which
 are defined over number fields instead of non-trivial division algebras.  Up to this point, the effort to construct an example of optimal canonical degree in the form of a lattice of first type has not been successful.

\ms
\ni{\bf 1.2.}  The idea of proof Theorem 1 is as follows.  The key observation is to relate a well-chosen fake projective plane to possible existence of a surface of optimal canonical degree.
 An appropriate normal cover of a fake projective plane of degree four gives the Euler number expected for a candidate surface.
We need to guarantee the vanishing of the first Betti number to achieve the correct dimension of the space of the canonical sections.  After this,
the main part of argument is to ensure that
the canonical map is generically finite and base point free, which turns out to be subtle.
  In this paper, we choose an appropriate covering corresponding to a congruence subgroup of the lattice associated
to an appropriate fake projective plane, which ensures that the first Betti number is trivial and the Picard number is one.  The latter condition makes the surface geometrically
simple for our arguments.  We divide the proof into three steps, proving that the rational canonical map is generically finite, that the map has no codimension one
base locus, and that the map has no codimension two base locus.  We make extensive use of the finite group actions given by the covering group.
In this process we have to utilize the geometric properties of the fake projective plane and relate to 
 finite group actions on a projective plane and on a rational line.   We also need to utilize vanishing properties in  [LY] of sections of
 certain line bundles which are numerically small rational multiples of the canonical line bundle, related to a conjecture on existence of exceptional objects in [GKMS].

 We would like to explain that the software package Magma was used in this paper, but only very elementary commands are used.  Starting with the presentation of our fake projective plane
 given in [CS], only one-phrase commands as used in calculators are needed  (see the details in \S3).
%The descriptions of the fundamental group of the surface was given in
%terms of generators and relations derived from the corresponding one of the fake projective plane given by [CS], see also the weblink there.  Simple
%commands from the software package Magma are used to derive information about the abelianization of the discrete groups involved.   Apart from these readily obtained information, the proof involved only
%classical algebraic geometric and group theoretical arguments.

More examples and classification of surfaces of optimal canonical degree arising from fake projective planes would be discussed in a forthcoming work with Ching-Jui Lai.

\ms
\ni{\bf 1.3.}  The author is indebted to Donald Cartwright for his numerous helps related to the use of Magma and lattices associated to fake projective planes.  He would also like to thank Ching-Jui Lai for many discussions,
comments and help on the paper, to Carlos Rito for spotting some errors in earlier drafts and making helpful comments,
and to Rong Du for bringing the problem to
his attention.  Part of the work was done while the author is visiting the Institute of Mathematics at the University of Hong Kong, and the author would like to express his gratitude for the hospitality of the
institute.

\bs
\begin{center}
{\bf 2. Preliminaries} 
\end{center}

\ms
\ni{\bf 2.1.}  For completeness of presentation, let us explain why the maximal degree is bounded from above by $36$ as is observed in [B].
Let $S$ be the rational image of $\Phi_{K_M}$.
Denote by $F$ and $P$ the fixed and movable parts of the canonical divisor $K_M$ respectively
and $\pi:\hM\rightarrow M$ the resolution of $P$.  Let $\pi^*P=F_{\hM}+P_{\hM}$ be the similar
decomposition on $\hM$.
Let $h^1(M)=\dim_{\bC}H^1(M,\cO_M), p_g=\dim_{\bC}H^2(M,\cO_M)$ and $\chi(\cO_M)$ be the
the arithmetic genus of $M$
respectively.  Then 
\begin{eqnarray*}
(\deg S) d =P_{\hM}^2&\leqslant& P^2\\
&\leqslant& K_M^2\\
&\leqslant& 9\chi(\cO_M)\\
&=&9(p_g-h^1+1)\\
&\leqslant &9(p_g+1).
\end{eqnarray*}
where the first two inequalities were explained in [B].
Hence 
$$d\leqslant 9(\frac{p_g+1}{\deg S}).$$
However, from Lemma 1.4 of [B], we know that $\deg S\geqslant p_g-2$.  We conclude that
$$d\leqslant 9(\frac{p_g+1}{p_g-2})\leqslant 36,$$
since $p_g\geqslant 3$ from the fact that the canonical mapping is generically finite.

Tracing back the above argument, it follows that the equality holds only if the fixed part of $M$ is trivial, $p_g=3$ and $h^1=0$, and
that the Miyaoka-Yau inequality becomes an equality.   From the work of Aubin and Yau, see [Ya], the latter condition implies that it
is the quotient of a complex two ball by a lattice in $PU(2,1)$.  Moreover,  we see that the canonical mapping is base point free
by tracing through the argument above.

Note that a complex two ball quotient is infinitesimally rigid from the result of Calabli-Vesentini [CV].  Hence such a surface cannot
be constructed from complete intersections.

\bs
\begin{center}
{\bf 3. Description of the surface} 
\end{center}

\ms
\ni{\bf 3.1.}  Recall that a fake projective plane is a smooth compact complex surface with the same Betti 
numbers as $P_{\bC}^2$. This is a notion introduced by Mumford [M] who constructed 
the first example. All fake projective planes have recently been classified into twenty-eight 
non-empty classes by the work of Prasad-Yeung in [PY].  Together with the work of  Cartwright-Steger [CS], we know that
there are precisely $100$ fake 
projective planes among those $28$ classes. It 
is known that a fake projective plane is a smooth complex two ball quotient $B_{\bC}^2/\Pi$ for
some lattice $\Pi\subset PU(2,1)$, and has 
the smallest Euler number among smooth surfaces of general type.  We refer the readers to
[R\'e] and [Ye2] for surveys about fake projective planes.

The fake projective planes $X=B_{\bC}^2/\Pi$  are classified in the sense of lattices.
In the notations explained in [PY], lattices $\Pi$ are constructed as a subgroup of a lattice $\oGamma$
which determined a class of fake projective planes classified.

\ms
\ni{\bf 3.2.}  In this paper, we are going to consider the following specific fake projective plane.
In the formulation of [PY], 
the surface has the same defining number fields
as Mumford's fake projective plane as constructed in {\bf 5.7, 5.11} of [PY], corresponding to $a=7, p=2$ in the notation there.  In particular,
two different lattices $\Pi$ representing fake projective planes with automorphism group of order $21$ are constructed.  Each such
lattice $\Pi$ is a congruence subgroup as explained in {\bf 5.11} of [PY] and is different from the one of Mumford.   
According to [PY], the associated maximal arithmetic lattice $\oGamma$ is an arithmetic lattice of second type in the sense that it is an arithmetic
lattice defined from
a non-trivial division algebra $\cD$ with an involution of second type $\iota$.  
%The fake projective plane that we consider is one such particular surface,
%denoted by $(a=7, p=2, \emptyset, D_3, 2_7)$ in the notation of [CS], the file `registry of surface' in the weblink
%\verb'http://www.maths.usyd.edu.au/u/donaldc/fakeprojectiveplanes/registerofgps.txt'. 

%The presentation of the lattice $\oGamma$ and $\Pi$ is given in [CS].  In this way, it was found in [CS] with Magma that
%$H_1(X)=\bZ_2^4.$   Hence there are normal subgroups $\Delta$ of $\Pi$ of order $4$ given by the kernel of the homomorphism
%$$\a:\Pi\rightarrow H_1=\bZ_2^4\rightarrow \bZ_2^2,$$
%where the last homomorphism is given by projections into any two of the four factors of $\bZ_2$.  Then $M=B_{\bC}^2/\Delta$ is a normal
%covering of $X$ of degree $4$.  However, for our later argument, we choose a concrete normal subgroup $\Sigma$ of $\Pi$ as follows.

\ms
\ni{\bf 3.3.}  The maximal arithmetic group $\oGamma$ to be used in this article corresponding to the class chosen above in [PY] (cf. Theorem 4.2, {\bf 5.9}, {\bf 5.11}).
A presentation of the lattice is found with a procedure explained by Cartwright and Steger in [CS] and details given by the 
file \verb'a7p2N/gp7 2generators reducesyntaxtxt' in the weblink of [CS],
 with generators and relations given by
\begin{eqnarray*}
\oGamma:=\langle z,b &|&z^7, (b^{2}z^{-1})^3, (bz^{-1}b^3z^2)^3,
(b^3z^{-2}bz^{-2})^3,
b^3z^{-2}b^{-1}z^2b^{-2}z,\\ && b^3z^3bz^2b^{-1}z^{-1}b^3z,zb^2z^{-2}b^{-1}z^{-1}b^{-3}zb^{-1}z^{-1}b^3z, \\&&bzb^5z^{-2}b^2z^2b^2z^{-2}b^2z^3
\rangle.
   \end{eqnarray*}
%We remark that the presentation on the right hand side are found by Cartwright and Steger by checking 

%\ms
%\ni(a) that it generates a subgroup of the maximal arithmetic group $\oGamma$ determined in [PY],
%and \\
%(b) that it has covolume $\frac1{12}$, 
 %which is precisely the covolume given by Prasad's Volume formula in [Pr] for $\oGamma$, cf. [PY] Proposition 3.5.
  
 % \ms
%  (a) and (b) imply that our group with presentation on the right hand side of the above expression generates $\oGamma$.  For (b), we observe that
%  the covolume is theoretically a rational number with denominator $12$ while the volume of the fundamental domain of an explicitly presented
%  group can be checked up to very high order of accuracy, as established by Cartwright and Steger in [CS].
  
The lattice associated to the fake projective plane is denoted by $\Pi$ and is generated by the subgroup of index $21$ in $\oGamma$ with generators given by
$$b^3, z^2bz^{-1}b^{-1}, (zbz^{-1})^3, zbz^{-1}b^{-1}z, zb^{-1}z^{-2}b, (bz^{-1})^3,$$
which is one of the candidates found by command {\tt LowIndexSubgroups} in Magma and is the one we used, denoted by $(a=7, p=2, \emptyset, D_3 2_7)$
in the notation of Cartwright-Steger (see file \verb'registerofgps.txt' in the weblink of [CS]), the first two entries correspond to $a=7, p=2$ in the number fields studied in [PY].  Denote by $X$ the resulting fake projective plane.
It follows that $H_1(X,\bZ)=\bZ_2^4$, which follows after applying the Magma
 command {\tt AbelianQuotient} to the presentation above.

Denote by $g_1,\dots,g_6$ the elements listed above.  
Magma  command {\tt LowIndexSubgroups} allow us to find a normal subgroup $\Sigma$ of index $4$ in $\Pi$ with generators given by
$$g_4
    ,g_5 g_1^{-1}
    ,g_6 g_2^{-1}
    ,g_1^{-2}
    ,g_2^{-2}
    ,g_3^{-2}
    ,g_5^{-1}  g_1^{-1}
    ,g_6^{-1} g_2^{-1}
    ,g_1 g_2 g_3^{-1}
    ,g_1  g_3 g_2^{-1}.
$$
The corresponding ball quotient is denoted by $M=B_{\bC}^2/\Sigma$.  In the next few sections, we would show that $M$ 
 is a surface with maximal canonical degree.  
 
\bs
\begin{center}
{\bf 4. Some  geometric properties of the surface} 
\end{center}

\ms
\ni{\bf 4.1.}  We collect some general information about the surface $M$.
\begin{lemm}
The ball quotient $M$ is a smooth unramified covering of degree $4$ of the fake projective plane $X$ satisfying the following properties.\\
(a). $b_1(M)=0$ and $H_1(M,\bZ)\cong \bZ_2^5\times \bZ_4$.\\
(b). Picard number $\rho(M)=1.$ \\
(c).  The lattice $\Sigma$ is a congruence subgroup of $\Pi$.\\
(d). The automorphism group of $M$ has order $\Aut(M)=A_4$, the alternating group of $4$ elements.\\
(e). $\Sigma\lhd\Pi$, $\Pi\lhd N_\Sigma$, $\Pi\lhd \oGamma$ with $|\Pi/\Sigma|=4$, $|N_\Sigma/\Sigma|=12$, and $|N_\Sigma/\Pi|=3,$
where $N_\Sigma$ is the normalizer of $\Sigma$ in $\oGamma$.\\
(f).  The action of $\bZ_3$ on $M$ descends to an action of $\bZ_3$ on $X$.\\
(g). The sequence of normal coverings $B_{\bC}^2/\Sigma\stackrel{p}\rightarrow B_{\bC}^2/\Pi\stackrel{q}\rightarrow B_{\bC}^2/\oGamma$
corresponds to normal subgroups $\Sigma\lhd\Pi\lhd\oGamma$, with covering groups $\Pi/\Sigma=\bZ_2\times\bZ_2$ and $\oGamma/\Pi=\bZ_7:\bZ_3$, the unique
non-abelian group of order $21$. 
%In fact, $q\circ p:B_{\bC}^2/\Sigma\rightarrow B_{\bC}^2/\oGamma$ is a normal covering.

\end{lemm}

\ni{\bf Proof} From the presentation of $\Sigma$ and Magma command {\tt AbelianQuotient}, we conclude that
$H_1(M)\cong \bZ_2^5\times \bZ_4$.  Hence (a) follows.

To prove (c), we consider the division algebra $\cD$ associated to our fake projective plane mentioned above.  Let 
$$V=\{\xi\in \cD:\iota(\xi)=\xi, \Tr(\xi)=0\}.$$
$V$ forms a vector space of dimension $8$ over $\bQ$.  $\oGamma$ has a representation on $V$, acting by conjugations.  Hence there is a natural homomorphism
$f:\oGamma\to SL(8,\bZ)$.
Considering reduction modulo $2$, 
there exists a homomorphism $f_2:\oGamma\to SL(8,\bZ_2)$ for which $|f_2(\oGamma)|=64\times 21.$
From Magma, we can check that the image $f_2(\Pi)$ of $\Pi$ has order 64, and so has index 21 in the image of $\oGamma$.
Recall that $\Pi$ has index 21 in $\oGamma$.  Hence  $\Pi$ contains the kernel of $f_2$ and is a congruence subgroup of $\oGamma$.
The author is indebted to Donald Cartwright for explaining the above procedure checking congruence property.

Consider a normal subgroup $\Sigma$ of $\Pi$ with index $4$ given by choice in the last section below.   From Magma again, the order of $f_2(\Sigma)$
is $16$ and hence is of index $4$ in $f_2(\Pi)$.
Again, as $\ker(f_2)\subset \Sigma$, we conclude that $N$ is a congruence subgroup. Hence (c) is true.  

Once we know that (c) is true, the facts about Picard number in (b) and $b_1(M)=0$ in (a) also follow from the work of Rogawski [Ro] and Blasius-Rogawski [BR], see also [Re].

For (d) and (e), we check by magma that the normalizer $N_\Sigma$ of $\Sigma$ in $\bar\Gamma$ is a subgroup of index $7$.  Hence we know that the automorphism group of $M$ given by
$N_\Sigma/\Sigma$ is a group of order $12$.  In fact, this corresponds to the group 
$(a=7, p=2, \emptyset, 2_7)$
in the notation of Cartwright-Steger in file \verb'registerofgps.txt' in the weblink of [CS], since that is the only group of right order in $\bar\Gamma$ supporting a
unramified covering of index $12$.   From Magma,  we check that the quotient group $H:=N_{\Sigma}/\Sigma$ is a non-abelian with $[H,H]=\bZ_3$
and actually $H= A_4$ after comparing with the library of small groups in Magma.  Magma also allows us to show that  $\Pi\lhd N_{\Sigma}$.
%Since $\Aut(M)$ is abelian, clearly the action of $\bZ_3$ on $M$ is the lift of some action of $\bZ_3$ on $X$.

(f) follows from the fact that $C=\bZ_2\times\bZ_2$ is a normal subgroup of $A_4$.   Recall that $\Sigma$ is a normal subgroup of $\Pi$ with quotient $C$ so that
we may write $\Pi=C\Sigma$.  Let $x\in B_{\bC}^2$.  By definition, for $\gamma\in \bZ_3<A_4,$ the action of
$\gamma$ at the $\Sigma$ cosets satisfies
$$\gamma(\Sigma x)=\gamma \Sigma\gamma^{-1}\cdot\gamma x=\Sigma(\gamma x).$$
We need to show the same is true for a $\Pi$ coset.  This follows from
$$\gamma(\Pi x)=\gamma(C\Sigma x)=\gamma C\gamma^{-1}\cdot \gamma \Sigma x=C\Sigma(\gamma x)=\Pi x$$
where we used the fact that $C$ is a normal subgroup of $A_4$.

(g) follows from the above description as well.

\qed

\bs
\ni{\bf Remark} As a consequence of the Universal Coefficient Theorem, the torsion part of the N\'eron-Severi group corresponds to the part in $H_1(M,\bZ),$
namely, $\bZ_2^5\times \bZ_4$.

\ms
\ni{\bf 4.2.}  We also recall the following result which is related to a conjecture of Galkin-Katzarkov-Mellit-Shinder in [GKMS].

\begin{lemm}
Let $H$ be the ample line bundle on $X$ on the fake projective plane $X$ as studied above, so that $K_X=3H$ as defined in [PY], {\bf 10.2, 10.3}.  Then \\
(a). $H^0(X,2H)=0$.\\
(b). There is no $\Aut(X)$ invariant sections in $H^0(X,2H+e)$, where $e$ is any torsion line bundle on $X$.
\end{lemm}

\ni {\bf Proof}   Part (a) follows from Theorem 1.3 or Lemma 4.2 of
Galkin-Katzarkov-Mellit-Shinder [GKMS], Theorem 1 of Lai-Yeung [LY], or Theorem 0.1 of Keum [K].   Part (b) follows directly from the proof of Theorem 1 of [LY].

\qed

\ms
\ni{\bf 4.3.}  Recall that from construction in \S3, $M$ is an unramified covering of a specific fake projective plane $X$ of index $4$.
Since $X$ is a fake projective plane, the Betti numbers and Hodge numbers of $X$ are 
the same as the corresponding ones on $P_{\bC}^2$.  It follows that $\chi(\cO_X)=1$.  Hence $\chi(\cO_M)=4\chi(\cO_X)=4$.
Since $h^1(M)=0$, it follows that $p_g=4-h^1(M)-h^0(M)=3$.  We conclude that $h^0(M,K_M)=3$.
Let $\{s_1, s_2, s_3\}$ be a basis of $H^0(M,K_M)$.
The linear system associated to the basis gives rise to a rational mapping
$$\Phi:x\dashrightarrow [s_1(x),s_2(x),s_3(x)]\in P_{\bC}^2.$$
Let $ S$ be the rational image of $\Phi$.  We know that $\dim_{\bC} S=1$ or $2$.

%We observe the following.

\begin{lemm} The action of the covering group $\cG\cong \bZ_2\times \bZ_2$ on $M$ induces an action of $\cG$ on $\Phi(M)\subset P_{\bC}^2$ through
the canonical rational map $\Phi:M\dashrightarrow \Phi(M)\subset P_{\bC}^2$.
\end{lemm}

\ni{\bf Proof}  From factorization of rational maps in complex dimension two, there exists a morphism $\pi:\hM\rightarrow M$, where $\hM$ is a sequence of blow-ups of $M$, and holomorphic map $f:\hM\rightarrow P_{\bC}^2$ such
that $f=\Phi\circ\pi$.

We observe that the action of $\cG$ on $M$ lifts to an action on $\hM$.  To see this, observe that the base locus of the canonical map
is invariant under $\Aut(M)$.   $\hM$ is obtained from $M$ from a series of blow-ups and we know that the induced action is biholomorphic outside of the
blown-up locus.  Since the transformation $\gamma\in \cG$ comes from a fractional linear transformation of an element in $PU(2,1)$, 
the transformation is locally linear around any fixed point and in particular lifts to the blown up divisors.   Hence $\cG$ acts holomorphically on 
$\hM$.

We define an action of $\cG$ on $P_{\bC}^2$ as follows.  For $\gamma\in \cG$ and $z\in P_{\bC}^2$ satisfying
$z=f(x)$, define
$$\gamma z=f(\gamma(x)).$$
To see that it is well-defined,  suppose $f(x)=[g_1(x), g_2(x), g_3(x)]$ for some $g_i(x)$, $i=1,2,3$, corresponding to a basis of the linear system
associated to $K_M$.
Assume that $f(x)=f(y).$ 
 Then since $g_1, g_2, g_3$ form a basis for the space of sections in
$K_M$, we conclude that there exists a constant $k$ such that  $g_i(x)=kg_i(y)$ for all $1\leqslant j\leqslant 3$.  Hence
$g(x)=kg(y)$ for all $g\in \Gamma(M,K)$, from which we deduce that $\gamma^*g_i(x)=\gamma^*g_i(y).$  We conclude that
\begin{eqnarray*}
f(\gamma x)&=&[\gamma^*g_1(x), \gamma^*g_2(x), \gamma^*g_3(x)]\\
&=&[\gamma^*g_1(y), \gamma^*g_2(y), \gamma^*g_3(y)]\\
&=&
f(\gamma y),
\end{eqnarray*}
from which we conclude that $\cG$ acts on $\Phi(M)\subset P_{\bC}^2$.

\qed

\bs
\begin{center}
{\bf 5. Generically finiteness} 
\end{center}

\ms
\ni{\bf 5.1.}  The goal of this section is to show that the rational mapping $\Phi$ is dominant.  First we make the following observation.
Recall as in [PY] that $K_X=3H_X$ for some line bundle on $X$ which corresponds to
a $SU(2,1)$-equivariant hyperplane line bundle $H$ on $\wtM\cong B_{\bC}^2$.  In the following, the descends of $H$ to $M$ and $X$ would be
denoted by $H_M$ and $H_X$ respectively, or simply $H$ when there is no danger of confusion.  In particular, $H_M=p^*H_X$.  Hence $K_M=3p^*H_X$.
As $K_M^2=36$ and the Picard number of $M$ is $1$, there are the following two different cases to consider, 

\ms
\noindent {\it Case (A)}, $H_M=2L$, where $L$ is a generator of
the Neron-Severi group of $M$ modulo torsion, or \\
{\it Case (B)}, $H_M$ is a generator of the Neron-Severi group.

\ms
\ni{\bf 5.2.}
\begin{lemm}  
%Let $M$ be the surface defined in \S3.
The canonical map $\Phi$ of $M$ is generically finite.
\end{lemm}

\ni{\bf Proof}  
Assume that $\dim_{\bC} S=1$.  We claim that the rational image $C=\Phi(M)$ has genus $0$.  Assume on the contrary that $C$ has genus
at least $1$.
As mentioned earlier, there exists a morphism $\pi:\hM\rightarrow M$, and a holomorphic map $f:\hM\rightarrow P_{\bC}^2$ such
that $f=\Phi\circ\pi$.  By Hurwitz Formula, the blown up divisors are mapped to a point on $C$ and hence  actually $\Phi$ extends across any possible base
point set of $\Phi$ to give a holomorphic $\Phi:M\rightarrow P_{\bC}^2$.  As $M$ has Picard number $1$ from Lemma 1, this leads to a contradiction since the fibers are contracted.
Hence the Claim is valid.

In general, we may write $K_M=F+P$, where $F$ is the fixed part and $P$ is the mobile part.
In our case here, from the claim, it follows that $C=\Phi(M)$ is a rational curve.  
Since $\dim(\Phi(M))=1$, as mentioned in [B] {\bf 1.1}, page 123, we may write 
\begin{equation}
K_M\equiv F+2Q,
\end{equation}
 where $F$ is the fixed part of $K_M$,
$2Q$ is the mobile part of $K_M$ and $Q$ is an irreducible curve. 
Here we denote the numerical equivalence of two divisors $A$ and $B$ by $A\equiv B$.

Our next step to prove the claim that $F$ is trivial.  Assume on the contrary that $F$ ins non-trivial.

Consider first {\it Case (A)}.
If $F$ is non-trivial, it follows that $F\equiv bL$, where $b$ is even and hence $b\geqslant 2$, which in turn implies that $Q\equiv cL$, where $c=3-\frac b2\leqslant 2$
from the decomposition of $K_M$ above.

Consider first the case that $b=2$ so that $F\equiv 2L$. It follows that $H_M\equiv 2L$ since $\rho(M)=1$.  Hence $H$ is the same as $2L$ up to a torsion line bundle in $\bZ_2^5+\bZ_4$, from Lemma 1 and the Universal Coefficient Theorem.

Hence $F\equiv H_M$ on $M$.  As $F$ on $M$  by definition is invariant under $\Aut(M)$, it descends to $X$ to give an effective divisor $G$ on $X$.  It follows that
$G\equiv H_X$ on $X$.  As $H_1(X,\bZ)=\bZ_2^4$, $G=H_X+e_2$ for some two torsion line bundle on $X$ from the Universal Coefficient Theorem.  This implies that
$2H_X=2G$ is effective on $X$, contradicting Lemma 2a.

%It follows that $2F=K_M-Z\equiv 2H$ and we may take $F=H$ or $F=H+e_4$, where $e_4$ is the four torsion line bundle in $\bZ_2^5+\bZ_4$.  In the first case, as $K_X=3H$,
%we get $Z=H$.  However,     In the second case with $F=H+e_4$, again as $K_X=3H$, we conclude that $Z=K-2F=H-2e_4.$
%However, this implies that $2H=2Z$ is an effective divisor.  and reach a contradiction again. 

The only other possibility is that $b=4, c=1$.  In such case, we would have $Q\equiv L$.   Hence we may choose $F$ to be a generator of the N\'eron-Severi group 
modulo torsion on $M$.  In such case, we may write $H=2Q+e$, where $e$ is a torsion line bundle corresponding to an element in $H_1(M,\bZ)=\bZ_2^4\times\bZ_4$ from Universal Coefficient
Theorem.  It follows that $K_M=3H_M=6Q+3e$.  Since $K_M=F+2Q,$ we conclude that $F=4Q+3e$.  Hence $F=2H_M+e$.  

As the canonical line bundle $K_M$ is invariant under the automorphism group of $M$, we know that the dimension one component $F$ of the canonical line bundle
is invariant under $\Aut(M)$.  It follows that $F$ descends as an effective divisor $G$ on the fake projective plane $X$.
The line bundle $H$ is clearly invariant as a holomorphic line bundle under $\Aut(M)$ from construction.  
It follows from $e=F-2H$ that $e$ is invariant as a holomorphic line bundle under $\Aut(M)$.  We conclude that $H^0(X,2H_X+e)\neq 0$ on $X$, 
since it contains the effective divisor $G$, where $p:M\rightarrow X$ is the covering map.
  Recall that from our setting, the coverings $B_{\bC}^2/\Sigma\rightarrow B_{\bC}^2/\Pi\rightarrow B_{\bC}^2/\oGamma$
corresponds to normal subgroups $\Sigma\lhd\Pi\lhd\oGamma.$  Hence from construction $G$ is invariant under $\Aut(X)=\oGamma/\Pi$.  This contradicts 
Lemma 2b.  Hence $F$ is trivial for {\it Case (A)}.

\ms
Consider now {\it Case (B)}.   In such case, as $K_M=3H_M$, equation (1) implies that $F\equiv H_M$.  Again, as argued earlier in {\it Case (A)}, $F$ on $M$  descends to $X$ to give an effective divisor $G\equiv H_X$ on $X$.  Furthermore, $G=H_X+e_2$ for some two torsion line bundle on $M$ so that
$2H_X=2G$ is effective on $X$, contradicting Lemma 2a.

Hence the claim about triviality of $F$ is proved.  We conclude that $K_M=P$.  In general, $P$ may have still have codimension two base point set, which is a finite number of
points in this case.   From equation (1), we may write $K_M=2Q$ for an effective divisor $Q$ on $M$, where $Q$ is the pull-back of $\cO(1)$ on $\Phi(M)\subset P_{\bC}^2$,
here we recall that $\Phi(M)$ is a rational curve as discussed earlier.  Now applying Lemma 3, we see that $\cG$ induces
an action on the rational image $\Phi(M)\subset P_{\bC}^2$.  As $\Phi(M)$ is a rational curve, from Lefschetz Fixed Point Theorem, $\cG$ has two fixed points
on $\Phi(M)$.  Let $a$ be such a fixed point on $\Phi(M)$.  The fiber $\pi(f^{-1}(a))$ above the fixed point $a$ corresponds to an effective divisor $Q_1$ in the class of $Q$ on $M$ as mentioned above.

Note that $K_M=2Q$ also implies that only Case (I) may occur, that is, $K_M\equiv 6L$, where $L$ is a generator of the N\'eron-Severi group on $M$ and hence that $Q_1\equiv 3L$.  On  the other hand, $Q_1$ as constructed is fixed by
$\cG$ as a set.  Hence $Q_1$ as a variety is invariant under the action of the Galois group $\cG$
and descends to $X$ to give rise to an effective divisor $R_1$ on $X$.  Note that $Q_1$ contains all base points of $K_M$ and hence the orbits of any base
point, which is assumed to be non-trivial.  Hence $Q_1=p^*R_1$ and is connected.  On the other hand, $R_1\equiv cH_X$ on $X$, where $1\leqslant c\leqslant 3$ is a positive integer.
Hence 
$$Q_1=p^*R_1\equiv cH_M\equiv 2cL,$$
which contradicts the earlier conclusion that $Q_1\equiv 3L$.  Hence $P$ has no base point set.

If follows that $\Phi$ is a morphism and fibers over a rational curve.  However, this contradicts the fact that $M$ has Picard number $1$.

In conclusion, $\dim_{\bC} S\neq 1$ and hence has to be $2$.
 
\qed

\bs
\begin{center}
{\bf 6. Codimension one component of base locus} 
\end{center}

\ms
\ni{\bf 6.1.}  The goal of this section is to show that there is no fixed component in the linear system associated to $K_M$.
\begin{lemm}
The base locus of $\Phi_{K_M}$ does not contain dimension one component.
%Moreover, the canonical mapping $K_M$ of $M$ has degree $36$.
\end{lemm}
\ni{\bf Proof} Let $L$ be the generator of the N\'eron-Severi group modulo torsion.   Since the Picard number is $1$,  we know that
$L\cdot L=1$ from Poincar\'e Duality.  Replacing $L$ by $-L$ if necessary, we may assume that $L$ is ample.  Now
we may write 
$$K_M=F+P,$$
where $F$ is the fixed part and $P$ is the moving part.  

We claim that $F$ is trivial.  Assume on the contrary that $F$ is non-trivial.

%Consider first {\it Case (A)} that 
%$K_M\equiv 6L$.  Recall as in [PY] {\bf 10.2-3} or from the last section that $K_X=3H$ for some line bundle on $X$ which corresponds to
%a $SU(2,1)$-equivariant hyperplane line bundle on $H$ on $\wtM\cong B_{\bC}^2$.

%Hence we have $F\equiv aL$ and $P\equiv (6-a)L$, where $1\leqslant a\leqslant 6$.  In particular, we may consider $F$
%as a global section of $\Gamma(M, aL)$.
%Since the base point set of $P$ has codimension at least $2$, and sections in $F$ does not move, we conclude that
%$a<6-a.$  Hence $1\leqslant a\leqslant 2$. 

From construction, the covering $p:M\rightarrow X$ is a normal covering of order $4$ and we may write $X=M/\cG$, where $\cG=\bZ_2\times\bZ_2$ is
a order $4$ group corresponding to deck transformation of the covering.  Hence $\cG$ is a subgroup of the automorphism
group of $M$.  From definition, $F$ is invariant under the automorphism group of $M$ and hence is invariant under $\cG$.
It follows that $F$ descends to an effective divisor $G$ on $X$.  As $X$ has Picard number $1$, we know that $K_X\equiv \beta G$ for some positive
rational number $\beta$, observing that $p^*K_X=K_M$ is numerically an integral multiple of $F$.  
%As $K_X=3H_X$, where $H_X$ is the generator of the N\'eron-Severi group modulo torsion on $X$, we conclude that $bG\equiv 3H_X$.   Hence  $b=1$  or $3$.  

From the remark in Section {\bf 4} and the descriptions in Section {\bf 3}, we know that the set of torsion line bundles on $X$ is given by $\bZ_2^4$.  Hence from $K_X=3H_X$, either 

\ms
\ni(I) $\beta G=3H_X$ and $K_X=\beta G$, or\\
(II) $\beta G=3H_X+e_X$ and $K_X=\beta G+e_X$, where $e_X$ is a $2$-torsion line bundle on $X$.

\ms
Case (I) cannot occur, since in such case there is a non-trivial section for $\Gamma(X,K_X)$, contradicting that $X$ is a fake projective plane.

Hence we only need to consider Case (II).  In such case, there are the following three subcases.

\ms\ni
(IIa) $G\equiv H_X$, or \\
(IIb) $G\equiv 2H_X$, or\\
(IIc) $G\equiv 3H_X$. 

\ms
In Case (IIa), $G=H_X+e_X$.  Hence $2H=2G$ is effective.  This is impossible from Lemma 2.

In Case (IIb), again from Lemma 2, we can rule out $G=2H$ and conclude that $G=2H_X+e_X$ for some two torsion line bundle $e_X$.  The argument of the last two paragraphs of \S5 leads to a contradiction.

For Case (IIc), we have $G=3H+e_X$, $K_X=G+e_X$.  There are a few subclasses.

\ms
\ni Case (IIci), $p^*G=F$ is irreducible.  In such case, $K_M=F+e_{M2}$, where $e_{M2}$ is a two torsion line bundle on $M$.   However,
as $K_M=F+P$, it follows that the movable part of $K_M$ is numerically trivial.  This is a contradiction.\\

\ni Case (IIcii), $p^*G=F_1+F_2$ consists of two irreducible components.  In such case, $F_2=\sigma F_1$ for some $\sigma\in\cG=\bZ_2\times \bZ_2.$
By taking dot product with a generator of the N\'eron-Severi group modulo torsion, we conclude that $F_2\equiv F_1$ and hence 
$F_2=F_1+e^\prime_{M2}$ for some two torsion line bundle $e^\prime_{M2}$.  In such case
$K_M=F_1+F_2+p^*(e_X).$  From construction, we know that $F_1+F_2\equiv 3H_M\equiv K_M$ on $M$.  Again, this leads to 
a contradiction since $P$ would then be numerically trivial.\\

\ni Case (IIciii), $p^*G=F_1+F_2+F_3+F_4$ consists of four irreducible components.  In such case, we can reach similar contradiction
by similar argument as above.  Alternatively, we see from similar argument as in the last paragraph that
$$K_M\equiv F_1+F_2+F_3+F_4\equiv 4F_1.$$
This leads to a contradiction since we either have {\it Case (A)}, $K_M\equiv6L$, where $L$ is a generator of N\'eron-Severi group modulo torsion, or
{\it Case (B)}, $K_M\equiv 3H_M$ with $H_M$ being the generator of the Neron-Severi group of $M$ modulo torsion.

We conclude that the base locus of $\Phi_{K_M}$ has no codimension one components.  

\qed

%We remark that an alternative proof of the Lemma 5 would also be implied by the proof of Lemma 6 below.

\bs
\begin{center}
{\bf 7. Zero dimensional components of the base locus} 
\end{center}

\ms
\ni{\bf 7.1.}   From Lemma 3, the Galois group $\cG$ of the covering $p:M\rightarrow X$ induces an action on $P_{\bC}^2$.  
Let $S$ and $T$ be the order two automorphisms generated by the first and the second factor of $\bZ_2$ on $\cG$ respectively.
From the results of [HL] (see also [S], [W]), as homology class of $P_{\bC}^2$ corresponds to
  the canonical class on $M$ is invariant under $\Aut(M)$, we know that the fixed point set of each of $\{S,T,ST\}$ consists of a line and an isolated point,
so that the three points form vertices of a triangle and the three line segments form the sides of the triangle.  Denote the triangle by
$\Delta_{P_1P_2P_3}$.  Hence we may assume that
$S$ fixes the point $P_1$ and the line $\ell_1$ is the line through $P_2$ and $P_3.$  Similarly for $S$ and $T$.  The vertices are the fixed points of
$\cG$.

Since a line on $P_{\bC}^2$ is defined by a linear equation $a_1x_1+a_2x_2+a_3x_3=0$ on homogeneous coordinate
$[x_1,x_2,x_3]\in P_{\bC}^2$, it corresponds to the zero set of a holomorphic section $s\in \Gamma(M,K_M)$.  Hence the
pull back of $\ell_i$ on $M$, defined by $\pi(f^{-1}(\ell_i))$ is given by the zero set of $s_i\in \Gamma(M,K_M).$

 \begin{lemm}  There is no zero dimensional component in $\cap_{i=1}^3Z_{s_i}$ for $s_i$ as defined above. 
\end{lemm}

The rest of the section is devoted to the proof of the lemma, which we resort to counting of intersection numbers and group actions.  For this purpose,
we first make some observations.

From Stein Factorization, we may decompose $f=g\circ h$ into holomorphic maps, where $h:\hM\rightarrow S$ has connected fiber, $f:S\rightarrow P_{\bC}^2$ is finite and $S$ is a normal surface. The degree of $f$ is the same as the degree of $g$. Since $f$ is generically finite, we know that there can at most be a finite number of dimension one fibers for $h$ and hence for $f$.
 Suppose $C$ is a dimension one fiber of $f$. Let $\hs$ be a section in  $\Gamma(\hM,P_{\hM})$. We make the following claim.
 
 \ms
 \ni{\it Claim}: $\hs\cdot C = 0$ and $\hs$ does not intersect $C$ if $\hs$ does not share a component with $C$.
 
 \ms
To prove the claim, we let $D$ be a hyperplane section on $P_{\bC}^2$ which avoids the set of points which are the image of all such contracted components $C$. From projection formula, $f^*D\cdot C = 0$. On the other hand, $f^*D\in \Gamma(\hM,P_{\hM})$. Hence $\hs\cdot C = f^*D\cdot C =0$.  This implies that $\hs$ does not intersect $C$ if $\hs$ does not
 share a component with $C$.

In the following we are going to apply the {\it claim} several times.  In our situation, since the Picard number of $M$ and $X$ are both $1$,  $C$ would descend to a divisor $C_1$ of  $\Gamma(X, H +\epsilon)$ or  $(X, 2H + \epsilon)$ for some $Aut(X)$-invariant torsion line bundle $\epsilon$ in the fake projective space $X$, which does not exists from the vanishing results in [LY].

\ms
\ni{\bf 7.2.}  We may assume that  $\pi:\hM\rightarrow M$ is a resolution of $M$  invariant under $\Aut(M)$, so that $f:\hM\rightarrow P_{\bC}^2$ is a morphism.
From construction
$\pi(f^{-1}(\ell_i))$, $i=1,2,3$, is invariant under $\bZ_2\times\bZ_2$ and hence is a $\bZ_2\times\bZ_2$-invariant section $s_i$ of $\Gamma(M,K_M)$.  Note that
they are linear independent by construction and hence span $\Gamma(M,K_M)$, which has dimension $3$.  Since each of
them is invariant under the Galois transformation group $\bZ_2\times\bZ_2$ of $p:M\rightarrow X$, each descends to a global section $t_i$ of $K_X+\tau$, where $\tau$ is a torsion line bundle.  Since
$p^*\tau=0$, we know that $\tau$ is a $\bZ_2\times\bZ_2$-torsion line bundle.  If $\tau=0$, we reach a contradiction since $H^0(X,K_X)=0$.  Hence
we conclude that $s_i\in \Gamma(X,K_X+\tau_i)$, where $\tau_i$ are non-trivial $2$-torsion line bundles.  We note that they span all the possible sections of
bundles of form $K_X+\tau$ in the orbit of $\bZ_7$ of sections of $\Gamma(X,K_X+\tau_i)$, where $\tau$ is a $2$-torsion line bundle, since $\dim(\Gamma(M,K_M))=3$.  
Here as mentioned in {\bf 3.2},  we know from the computation of Cartwright and Steger that $H_1(X,\bZ)=\bZ_2^4$, hence the bundle $K+\tau_i$ is
invariant under $\bZ_7$ as a line bundle.  Now for sections of $\Gamma(X,K_X+\tau_i)$, if a section is not invariant, the space would have dimension greater than $1$, which when lifted
to $X$ and taken together with $s_1,s_2,s_3$, would lead to $\dim(\Gamma(M,K_M))>3$.

  Let $B=p(A)$.  Since $\{s_1, s_2, s_3\}$
is $\bZ_2\times \bZ_2$ invariant, the zeros divisors $t_i$ all pass through each point of $B$ on $X$.  As $K_X\cdot K_X=9$, it follows that 
$B$ has at most $9$ points.

In the following we would denote by $\hs_i$ the proper transform of $s_i$ in the $\Aut(M)$-invariant minimal resolution $\hM$ of $M$ associated to the birational map $\Phi_{K_M}$.
The covering  map $p$ induces an isomorphism of a small neighborhood of base point of $s_i, i=1,\dots, 3$ to a small neighborhood of $t_i, i=1,\dots,3$.  
For convenience, we would denote by $\hht_i$ the proper transform of $t_i$ on $\widehat{X}$, the induced modification of $X$ corresponding to $\hM\rightarrow M$.

\ms
\ni{\bf 7.3} In terms of the notation of {\bf 2.1},
we note that $P^2=p^*P\cdot P_{\hM}=F_{\hM}\cdot P_{\hM}+P_{\hM}^2.$
The sequence of estimates of degrees can be written as 
\begin{equation}
\deg \Phi=\deg (f)=P_{\hM}^2=P^2-F_{\hM}\cdot P_{\hM}\leqslant P^2\leqslant K_M^2=36.
\end{equation}

\ms
\ni{\bf 7.4} Now we recall that $\Aut(X)$ is the abelian group $G=7:3=\bZ_7\rtimes\bZ_3$ of order $21$, where $\bZ_3$ acts on $\bZ_7$ by
a homomorphism $\bZ_3\rightarrow \Aut(\bZ_7)$.  $G$ has a normal Sylow subgroup of order $7$, denoted by $\bZ_7$.  There are also seven Sylow
subgroups of order $3$.  $\bZ_7$ has three fixed points on $X$, and each Sylow $3$-subgroup has $3$ fixed points on $X$, according to a result of Keum and
Cartwright-Steger.
Let $1\neq\gamma\in G$.  Note that $\gamma^* t_i$ would be another section of some $K_X+\tau$, where $\tau$ is a $2$-torsion.
As mentioned in the last paragraph, from dimension considering, it follows that $\tau$ has to be one of the $\tau_i, i=1,2,3$ mentioned
earlier, and $\gamma^*t_i$ has to pass through each point of $T$ as well.  As $|G|=21$ and the set $B$, which has cardinality at most $9$,
is invariant under an automorphism of $M$,
we conclude that each point $Q$ in $B$ is actually fixed by some element $\gamma\in G$.  In our case, there is a unique subgroup $\bZ_7$ of order $7$ and seven
subgroups $\bZ_3$ of order $3$ acting on $X$.   We consider the subgroup $\bZ_3$ of $\Aut(X)$ descended from $\Aut(M)$ as mentioned in Lemma 1(f).  The group of
order $3$ and the group of order $7$ generates $\Aut(M)$.
There are two cases to consider, 

\ms\ni
Case I: $Q$ is a fixed point of a subgroup $H$ of $\Aut(M)$ isomorphic to $\bZ_3$, and \\
Case II: $Q$ is a fixed point of the subgroup of $\Aut(M)$ isomorphic to 
$\bZ_7$.

\ms
\ni{\bf 7.5} Consider first Case I.  For simplicity, we just call the group involved $\bZ_3$.
 Assume now that a point $Q_1$ in $B$ is a fixed point of $\bZ_3$.    Then $Q_1$ lies on $t_1$ and is not fixed by $\bZ_7$, as it is well-known that
no point on
$M$ is fixed by the whole group $G=7:3$.  Hence the orbit of $Q_1$ has seven points $Q_i, i=1,\dots,7$ and all lies on $t_1$.  
As $t_1\cdot t_2=9$, we conclude that apart from the seven points $Q_1$ which are base locus, $t_1$ intersects $t_2$ either twice at a point or once at two points,
which we denote by $W$.  Since each point in $p^{-1}(W)$ is mapped to the point $\ell_1\cap\ell_2$ on $P_{\bC}^2$ and $\ell_1$ intersects $\ell_2$ in
simple normal crossing, we conclude that the degree $\deg(f)\geqslant 8$.  Here we recall from the discussion following the claim in {\bf 7.1}
that  $p^{-1}(W)$ does not
contain any one dimensional component.
Recall that each $t_i$ is fixed as a set by
$\bZ_7$.  Hence $A$ has $28$ points $R_i, i=1,\dots,28$ on $M$.  After resolving $A$, the base point set of $K_M$, each $R_i$ gives rise to an exceptional curve
$F_i$, which intersects each proper transform $\hs_i$ of $s_i$.  $\hs_1\cdot F_i>0$ for each $i=1,\dots,28$.  Hence $F\cdot P_{\hM}\geqslant 28$  in (2).  From (2)  
it follows that  $\deg(f)\leqslant 36-28=8$.   Hence we conclude that $\deg(f)=8$ and each proper 
transform of $\hs_i$ intersects $F_i$ only once for each $i$.  Moreover, the exceptional divisor over each $R_i$ is a single rational curve $F_i$ for each
$i=1,\dots,28$.

Now for each fixed $F_i$, $f(F_i)$ intersects $\ell+1$.  Let $x\in f(F_i)\cap \ell_1$.  
%Suppose that they meet at a point $x\in\ell_1$.  Note that $x$ cannot lie in the intersection of $\ell_1$ and $\ell_j$ for $j\neq i$, for otherwise
%the extra degree given by $F_i\cap \hs_j$ would lead to $\deg(f)>8$, contradicting the earlier conclusion.
We observe that $f^{-1}(x)$ contains $\gamma(F_i\cap f^{-1}(x))$ for all $\gamma\in \bZ_2\times \bZ_2$.
As the degree $\deg(f)=8$,  the degree of the curve $\gamma(F_i)$ in $P_{\bC}^2$ is either $1$ or $2$.

Consider the action of  a subgroup $\bZ_3$ of $G$ on $X$.  Either 

\ms\ni (a):  it leaves $t_i$ invariant as a set for all $i$, or\\
(b):  permutes among $t_i$, $i=1,2,3$.  

\ms
Consider first Case (a).  From [Su], [HL] or [W], as the homology class of $P_{\bC}^2$ corresponding to
  the canonical class on $M$ is invariant under $\Aut(M)$, we know that the fixed point set of $H\cong\bZ_3$ on $P_{\bC}^2$ consists either of a line and a point, or three points.  
We claim that the first case cannot happen.  Otherwise the line has to be one of $\ell_i, i=1,\dots, 3$ as the fixed point set contains $f(\hs_i\cap \hs_j)$ for $i\neq j$.  However,
if say $\ell_1$ is fixed by $\bZ_3$, it implies that $\hs_1$ is fixed pointwise under $\bZ_3$, since we know that the degree of $f$ is $8$, which is not divisible by $3$. 
Now we used the fact that $M$ as an unramified covering of $X$ is an arithmetic ball quotient division algebra and hence supports no totally geodesic curves, which in 
turn implies that a non-trivial finite group action on $M$ has only a finite number of fixed points, cf. [Ye1], p.19-21.  In particular, there is no fixed point of $\bZ_3$ on $\Phi^{-1}(y)$ for a generic
$y\in \ell_i\subset P_{\bC}^2$.
We conclude that
$\widehat{t}_1$ and hence $t_1$ is fixed pointwise  under $\bZ_3$, which is a contradiction since $\bZ_3$ has isolated fixed points on $M$.  Hence the claim is proved.

Hence the induced action of $\bZ_3$ on $P_{\bC}^2$ has three fixed points.  We aso know that on each rational line $\ell_i$, the induced action of $\bZ_3$ has two fixed points.  
Since there are three lines, it follows that the fixed points
of $\bZ_3$ has to be the three points $P_1, P_2, P_3$ corresponding to $\ell_i\cap\ell_{i+1}$.  

In terms of our earlier notation, we note that each $Q_i, i=1,\dots,7,$ which is a $\bZ_7$-orbit of $Q_1$, lies in $B$ as well, since  the divisors 
$t_j, j=1,2,3$ are invariant under $\bZ_7$ as we note earlier.  Note that the pull-back of each $t_j$ is just a single irreducible component $s_j$, which follows from 
Lefschetz Hyperplane Theorem.  In fact, any extra component would lead to $\dim(\Gamma(M,K_M))>3$ and a contradiction as well.  Hence each $R_j\in A$ for $j=1,\dots,28$.

Observe from Lemma 1 that $A_4$ is the automorphism group of $M$  and hence contains four $\bZ_3$-subgroups $H_i, i=1,\dots,4,$ permuted under conjugation by $\bZ_2\times \bZ_2$.  Moreover, the action of $H_i$ descends to $X$.  Let $H_1$ be the $\bZ_3$-group studied in the last paragraph.  On $X$, the points $g Q_i$, where $g\in \bZ_7$, are 
fixed by the $\bZ_3$-subgroup $gH_1g^{-1}$ of $\Aut(X)$.  Since $\Aut(X)=7:3$ contains precisely seven such subgroups under conjugacy of elements of $\bZ_7$, we know 
that 
four of the seven groups are given by $H_i, i=1,\dots,4$.  Consider now the four points $R_{1i}\in p^{-1}(Q_1),$ where $i=1,\dots 4.$  The set is invariant under $H_1$.  Hence we may
assume that $R_{11}$ is fixed under $H_1$.  We claim that each of $R_{1i}, i=2,3,4,$ is invariant under some $H_j$ for $j=2,3,4$.  This follows from the fact that the deck
transformation of the covering $p:M\rightarrow X$ is precisely $\bZ_2\times\bZ_2$, and that $hR_{11}$ is fixed by $hH_1h^{-1}$ for $h\in \bZ_2\times\bZ_2$.  This argument actually folds for
$R_{ji}, i=1,\dots,4$ for each $j=1,\dots,4$.  Hence there are sixteen such points $R_{ji}$.  Rename them as $R_i, i=1,\dots,16$.
It follows that $f(F_i\cap \hs_j)$ lies on $\ell_j$ for $i=1,\dots,16$ and $j=1,2$.

Recall from earlier discussions that the action of $H_k\cong \bZ_3, k=1,\dots,4$ on $P_{\bC}^3$ has three fixed points $P_k, k=1,2,3$.  From earlier discussions, we also know that $F_i\cap\hs_j$ only at one point.
Since both $F_i$ and $\hs_j$ are invariant under $\bZ_3$, it follows that $F_i\cap \hs_j$ for each $i$ and $j$ is invariant under $\bZ_3$.  Hence the same is
true for $f(F_i\cap \hs_j)$.  It follows that  for $j=1,2$ and $i=1,\dots,16$, $f(F_i\cap \hs_j)$ is one of the three fixed points mentioned earlier.  Since they also lie on
$\ell_1$ and $\ell_2$ by definition, it follows that $f(F_i\cap \hs_j)=\ell_1\cap\ell_2=P_3.$
Since there are at least $16$ points in the preimage of $P_3$ as constructed,  this contradicts $\deg(f)=8$ derived earlier.  Here we have used the {\it Claim} in \S7.1.

\bs
  Consider now Case (b).  In terms of earlier notation 
  $H_i\cong \bZ_3$ induced an action on $P_{\bC}^3$, the image of $\Phi$.  From construction, we know that $H_i$ leaves the three lines $\cup_{j=1}^3\ell_j$ invariant as a set, and permutes
  the three lines.
   From the results of [Su], [HL] or [W], $H_i$ acts as elements in $U(3)$ and the fixed point set consists either of  (i) three fixed points, or (ii) a point and a line $L$.  First we observe that (ii) cannot happen,
  for otherwise $L$ intersects $\ell_1$ and there is a fixed point of $H_i$ on $\ell_1$.  This implies that the fixed point has to be either $\ell_1\cap\ell_2$ or $\ell_1\cap \ell_3$.
  In the first case, $H_i$ has to permute between $\ell_1$ and $\ell_2$, which is not possible as $H_i$ has order $3$ and does not leave $\ell_1$ invariant.  Similar 
  contradiction arises in the second case.  Hence only (i) occurs.   Choose homogeneous coordinates on $P_{\bC}^2$ so that $\ell_1$ be defined by $Z_1=0$
  and $Z_2=\gamma Z_1, Z_3=\gamma^2 Z_1$, where $\gamma$ is a generator of $H_1$.  It follows that we may represent $\gamma$ in terms of our basis 
   $$
  \gamma=\left(\begin{array}{ccc}
  0&0&1\\
  1&0&0\\
  0&1&0
  \end{array}\right).$$
  For $H_i, i=2,3$, a generator $\gamma_i$ has to be of form 
  $$
  \gamma_i=\left(\begin{array}{ccc}
  0&0&\theta_{i1}\\
  \theta_{i2}&0&0\\
  0&\theta_{i3}&0
  \end{array}\right)$$
  where $\theta_{ij}, j=1,2,3$ are third roots of unity.  It follows from direct computation that $\theta_{i1}\cdot\theta_{i2}\cdot\theta_{i3}=1$ so that if we write
  $\theta_{ik}=\theta_{i1}\cdot\omega_{ik}'$ for $k=2,3$, we have $\omega_{i3}'=(\omega_{i2}^\prime)^2.$
  Hence in terms of the chosen homogeneous coordinates on $P_{\bC}^2$, the fixed points for each $H_i$ is given by 
  $U_1=[1,1,1], U_2=[1,\eta,\eta^2], U_3=[1,\eta^2, \eta]$, where $\eta$ is a third root of unity.

 Recall also that $H_i$ has a fixed point at $Q_i$ in our earlier notation for $i=1,\dots,4$.   Consider now $p^{-1}(Q_1)=R_{1j}, j=1,\dots 4$.
  We have assumed that $R_{11}$ is fixed by $H_1$.  Now since the exceptional divisor $F_{1}$ consists of a rational curve and we have 
  an action of $H_1\cong \bZ_3$ acting on $F_{1}$, there are at least two points on $F_{1}$ fixed by $H_1$.  Similarly, as the three points $R_{1j}, j=2,3,4$ are obtained
 from action of $\bZ_2\times \bZ_2$ on $R_{11}$,  we see that each $R_{1j}$ is fixed by some conjugate of $H_1$ and hence by one of $H_2, H_3, H_4$.  In other words,
 each of the four points $R_{1j}, j=1,\dots,4$ is fixed by precisely one $H_k$ for some $k=1,\dots,4.$  This holds for all the $16$ points $p^{-1}(Q_j)=\{R_{j1}, R_{j2},R_{j3},R_{j4}\}$, $j=1,\dots,4$.  Each of them gives rise to two fixed points of some $H_k$ on the exceptional divisor $F_{i}$.  Hence there are altogether $32$ such points.  Now 
the action of each of the four groups
 $H_i$, $i=1,\dots,4$ on $P_{\bC}^2$ has fixed point set given by $\{U_1, U_2,U_3\}$.
It follows that the degree of the mapping $\Phi$ is at least $32/3$, noting that there may be other points in the preimage.  Since $32/3>8$, this contradicts our earlier conclusion
that $\deg(\Phi)=8.$

\ms
\ni{\bf 7.6} Consider now Case II and again denote by $\bZ_7$ the unique $\bZ_7$ subgroup of $\Aut(M)$.  Assume now that a point $Q_1$ in $B$ is a fixed point of $\bZ_7$.    Under the action of $1\neq\gamma\in H_i, i=1,\dots 4$, a $\bZ_3$ subgroup 
 descends to $M$
mentioned earlier, we know that  that $\gamma Q_1\neq Q_1$ and hence has to be fixed by a conjugate of $\bZ_7$ subgroup of $G$.  As such a Sylow $7$-subgroup is unique,
the group is just the $\bZ_7$ group studied.  Hence $\gamma Q_1$ is fixed by $\bZ_7$ as well.  Moreover, the same argument implies that $\gamma Q_1$ lies in the base locus of
$t_j, j=1,2,3$ and hence $\gamma Q_1\in B$.  It follows that all
the three fixed
points of $\bZ_7$ on $X$ lie in $B$.
As discussed earlier, each $t_i$ is fixed as a set by $\bZ_7$.
Then $Q_1$ lies on $t_i$ and is not fixed by $\bZ_3$.  Hence
its orbit by $\bZ_3$ consists $Q_1$, $Q_2$ and $Q_3$ lying on $t_i$ for $i=1,2,3$.  This leads to $12$ base points $R_j$, $j=1,\dots,12,$ on $M$ after pulling back by $\pi$.  Resolving each 
base point in a $\Aut(M)$-invariant manner, it follows as before that $F\cdot \hs_i=F\cdot P_{\hM}$ is a positive multiple of $12$.    Note that for each of $i=1,2,3$, the behavior of $F\cdot s_i$ at all the points $R_j$ are all the same for all $1\leqslant j\leqslant 12$, since $s_i$ is invariant under each $H_j, j=1,\dots,4$.  Since $F_i\cdot P_{\hM}>0$ for each  
irreducible component $F_i$ of $F$, we conclude from (2) that $F\cdot P_{\hM}|_{R_j}=1$ or $2$, and $F|_{R_j}$ can have either 

\ms\ni
{\it Case (a)}, one component, or\\
{\it Case (b)}, two components.

Moreover,
\begin{equation}
\deg\Phi=\deg(f)=36-F\cdot P_{\hM}\leqslant 24.
\end{equation}

%We conclude from (1) as before that 
%$F\cdot P_{\hM}\geqslant 12$ and hence $\deg(f)\leqslant 36-12=24$. 
\ms
Now we observe that $t_1$ and $t_2$ cannot intersect at
any other points apart from base locus.  Otherwise there would be at least $7$ such points in the orbit of $\bZ_7$ on $t_1$.  This leads to 
$28$ points of the intersection of $s_1$ and $s_2$ on $M$.  Unless $\hs_1$ and $\hs_2$ share a component $C$ which is $f$ exceptional and mapped to the point $\ell_1\cap\ell_2$, the claim in \S7.1 implies that $\deg(f) > 28$, contradicting $\deg(f)\leqslant 24$. However if such a component $C$ exists, as $C$ does not have a component in the exception divisor of $\pi$
as studied in \S7.1, we conclude that $s_1$ and $s_2$ share some $C_1$ from $M$. This implies that $t_1$ and $t_2$ share some component $C_2$ on $X$. In such a case, $C_2$ is a section of  $(X, H + \epsilon)$ or  $(X, 2H + \epsilon$) for some $\Aut(X)$-invariant torsion line bundle $\epsilon$ in the fake projective space $X$, which does not exist from the vanishing results in [LY].

%Since $\bZ_7$ leaves each $t_i$ and hence $s_i$ invariant as a set, the same is true on $\ell_i$.  
%Observe now that $\bZ_7$ acting on $\ell_i$, which is biholomorphic to $P_{\bC}^1$,  has exact $2$ fixed points.  We also know that the action of $\bZ_7$ on 
%$P_{\bC}^2$ has fixed point  set consisting either of (i) three fixed points, or (ii) one fixed line and one fixed point, according to the results of [S], [HL] or [W].
%We claim that none of the $\ell_i$ is fixed pointwise  by $\bZ_7$.
%Assume on the contrary that $\ell_1$ is fixed point-wise by $\bZ_7$.  Note that $s_1$ cannot be fixed pointwise by $\bZ_7$, for otherwise $t_1$ would be totally geodesic
%(cf. [Y]), which is a contradiction since lattices associated to fake projective planes are arithmetic lattices of type two and does not support any totally geodesic curves.
%It follows that each generic points on $\ell_1$ corresponds to a multiple of $7$ points on $t_1$.  Since $\ell_1$ is fixed by a factor $\bZ_2$ of $\bZ_2\times \bZ_2$,
%we conclude that $\deg(f)$ is a multiple of $14$ and hence has to be $14$ from the last paragraph.  However, in such case (2) leads to $F_{\hM}\cdot P_{\hM}=22$.  
%On the other hand, as the $12$ points of $R_i$ are permuted by the automorphism group, and we see that $66$ is not divisible by $12$, we lead to a contradiction.
%The claim is proved.

%From the claim, it follows that the fixed points can only occur at the three points corresponding to $\ell_i\cap\ell_{i+1}$, namely one of the $P_k$
%for $k=1,2,3$.

Note that the three fixed points of $\bZ_7$ on $X$ are permuted by any subgroup isomorphic to $\bZ_3$ in $\Aut(X)$, and so does the base points of $K_M$ under the action
of $\Aut(M)$.  Hence the behavior of the base locus at the $12$ points of base locus on $X$ are the same.  Consider one such base point $R_a$.  Suppose
$F_{ai}, i=1,\dots,N$ are the irreducible components of the resolution $F_a$ of the point $R_a$ so that the proper transform of $\pi^*\Phi_{K_M}$ is base point free.  We note that
a resolution in a small neighborhood of a point $R_a\in A$ can be considered as the resolution of the corresponding point $Q_a\in B$, since the mapping $p:M\rightarrow X$ is
etale.  Hence by doing surgery, we may assume that there is a resolution $\pi:\hat X\rightarrow X$ for which an exceptional fiber $G_a$ at $Q_a$ is isomorphic to
an exceptional fiber $F_a$ at $R_a$.  Similarly, we let $\hht_j$ be the proper transforms of $t_i$.
Now since $Q_a$ is a fixed point of $\bZ_7$ on $X$, $\bZ_7$ acts on the exceptional fiber $G_a$ at $Q_a$.  
The induced
action of $\bZ_7$ should leave each $G_{ai}$ which intersect with some $\hht_j$ invariant.  Otherwise, there would be at least seven such components, giving
rise to $G_a\cdot\hht_i\geqslant 7.$   This is translated to the conclusion that $F_a\cdot\hs_i\geqslant 7$ on $\hM$.  Since there are $12$ such base points, it would lead to $F\cdot K_M\geqslant 12\cdot 7$ which violates (3).

\ms
Consider now {\it Case a}.   There is only one irreducible component in $G_a$ at $Q_a$.  Since $G_a$ is a rational curve, $\bZ_7$ has two fixed points only. 
Since $G_a\cap \hht_i$ is a fixed point, we may assume that the two fixed point are 
$\hQ_1=G_a\cap\hht_1$, and $\hQ_2=G_a\cap\hht_2=G_a\cap\hht_3$.   This is reflected correspondingly for $\hs_i$ on $R_j$.
Since $F_a\cdot\hs_i=F_a\cdot P_{\hM}$ for each $i$, this number can either be $1$ or $2$ from (3).  $G_a\cdot\hht_i$ cannot be $2$, for otherwise
the intersection of $\hs_2$ and $\hs_3$ at $F_a$ satisfies
$\hs_2\cdot\hs_3|_{F_a}=\hht_2\cdot\hht_3|_{G_a}\geqslant 2$, where the notation refers to intersection along $F_a$ or $G_a$.  Since there are twelve such points $R_a$, by looking at the preimage of 
$\ell_2\cap\ell_3$,
this implies that $\deg\Phi\geqslant 24$, contradicting (3) since $F\cdot P_{\hM}=F\cdot \hs_i=24$ in such case.  Hence we conclude that $F_a\cdot\hs_i=1$ for each $a=1,\dots,12$.
  In particular, we 
conclude from this and $F\cdot \hs_i=12$ that each $\hs_i$ intersects $F_a$ normally for each $i=1,\dots,3$ and $\hs_2$ intersects $\hs_3$ normally.
This implies that on $M$, $s_1$ intersects $s_2$ and $s_3$ transversally respectively, and $s_2$ intersects $s_3$ with multiplicity two at $R_a$.  Hence 
$t_1$ intersects $t_2$ and $t_3$ transversally respectively, and $t_1$ intersects $t_2$ with multiplicity two at $Q_a$.
This means that $(t_1\cdot t_2+t_2\cdot t_3+t_3\cdot t_1)|_{Q_a}=4,$ where $t_k\cdot t_l|_{Q_a}$ refers to multiplicity of  intersection of  $t_k$ and $t_l$ at $Q_a$.
Recall now that on $X$, the zero divisors $t_i\cdot t_j=K_M^2=9$
for $i\neq j$.  Hence
$$27=t_1\cdot t_2+t_2\cdot t_3+t_3\cdot t_1=\sum_{i=1}^3(t_1\cdot t_2+t_2\cdot t_3+t_3\cdot t_1)|_{Q_i}=12,$$
which is a contradiction.  

%Since all the intersections of $\hs_i$, $\hs_j$ for $i\neq j$ are mapped into one of the three intersection
%points $P_k$ of $\ell_\alpha$ on $P_{\bC}^2$, we conclude that the degree of $\Phi$ is at least $12\cdot 9/3=36$, which violates 
%$\deg(f)\leqslant 28$.  Here again we used the claim in \S7.1.

\ms
Consider now {\it Case b}.   In this case, an exceptional fiber $G_a$ at $Q_a$ consists of two irreducible components $G_{a1}$ and $G_{a2}$
meeting at a point $W_{a0}$ on $\hM$.  $W_{a0}$ is fixed by $\bZ_7$.  Denote by $W_{ai}$ the other fixed point of $\bZ_7$ on $G_{ai}$, $i=1,2$.
From (3) as before, we know that $\deg(\Phi)=12$ and $F\cdot P_{\hM}=24$.  As there are twelve points $R_i, i=1,\dots,12$ under consideration, we conclude that $F_a\cdot \hs_i=2$
for all $a=1,\dots,12$ and $i=1,2,3$.  As in {\it Case a}, $G_a$ meets $\hht_i$ only at one of the three fixed points of $\bZ_7$.  
%If none of them meet $W_{a0}$, the
%resolution is not minimal, since we may consider blowing down one component and the resulting $\Phi$ is still base-point free, which would be ruled out by 
%the discussions in {\it Case a}.  Similar argument applies if
If
$W_{a1}$ does not lie in at least one of $\hht_j, j=1,2,3$, 
%all the $\hht_i$ can only meet $G_a$ at the same component $G_{a1}$ or $G_{a2}$. However, 
as $P_{\hM}\cdot G_{a1}>0$, it follows that all $\hht_i, i=1,2,3$ intersects $G_{a1}$ at the point $W_{a0}$, which however contradicts that
$P_{\hM}$ is base point free.  Similarly, if $W_{a2}$ does not lie in one of $\hht_j, j=1,2,3$, it leads to the same contradiction.
If on the other hand $W_{a0}$ does not lie in at least one of $\hht_j, j=1,2,3$, all the $\hht_i, i=1,2,3$ meet $F_{a1}$ at the two points $W_{a1}$.
Again it follows that $W_{a1}$ is a base point of $P_{\hM}$ and leads to a contradiction.

Hence after renaming index if necessary,
we may assume that 
$W_{a1}\in G_a\cap \hht_1, W_{a2}\in G_a\cap \hht_2$ and $W_{a0}\in G_a\cap\hht_3$.
However, this implies correspondingly that $F\cdot\hs_3\geqslant\sum_{a=1}^{12}F_a\cdot \hs_3=24$.
From (3), it follows that $F\cdot\hs_3=24$.  Hence we conclude that $F\cdot \hs_i=F\cdot P_{\hM}=24$ for
$i=1,2$.
This implies that $\hht_i$ intersects $G_{ai}$ with multiplicity $2$ at $W_{ai }$ and hence $s_i$ intersects
$F_{ai}$ to multiplicity $2$, where $i = 1,2$. Now from the paragraph immediately after
(3), we conclude that $\hs_i$ cannot intersect $\hs_j$ at any point except for the union of the 
fibers $F$, which implies that $\hs_i\cdot \hs_j=0$ for $i\neq j$ from the discussion above.  This however
contradicts the fact that $\hs_i\cdot \hs_j=P_{\hM}\cdot P_{\hM}>0$.

In conclusion, both {\it Case (a)} and {\it Case (b)} leads to a contradiction.  Hence Case II does not occur.  This concludes the proof of Lemma 2.

\qed

\ms
\ni{\bf 7.7}   {\bf Proof of Lemma 6} 

From the discussions in {\bf 7.4}, every point in  the base locus has to be in the $\bZ_2\times \bZ_2$ orbit of
one of the fixed point set of either a subgroup of order $7$ or $3$
of the automorphism group of $\Aut(X)$. 
The discussions in {\bf 7.5} and {\bf 7.6} implies that there is no base locus corresponding to the fixed point set of $\Aut(M)$.
Lemma 6 follows.

\qed

\bs
\begin{center}
{\bf 8. Conclusion of proof} 
\end{center}

\ms
\ni{\bf 8.1.}  The discussions of the previous few sections can be summarized into the following proposition.
\begin{prop}
The linear system associated to $\Gamma(M,K_M)$ is base point free and the image of $\Phi_{K_M}$ is $P_{\bC}^2$.
\end{prop}

\ni{\bf Proof}  From Lemma 5, the base locus of $K_M$ is of dimension $0$.  From Lemma 6, we know that it is base point free.  From Lemma 4, 
we know that the image of $\Phi_{K_M}$ has complex dimension $2$ and hence has to be $P_{\bC}^2$.
%we are done if we know that a particular basis of
%sections of $\Gamma(M,K)$, which are fixed by some elements in the induced action of $\Aut(M)$,  have irreducible zero divisors.  
%From Lemma 6, Lemma 5 are readily applicable even if the sections 
%are reducible.

\qed

\ms
\ni{\bf 8.2.}  We can now complete the proof of our main result.

\ms  
\ni{\bf Proof of Theorem 1}  

We use the fake  projective plane $X=B_{\bC}^2/\Pi$ with $\Pi$ as given in Section {\bf 3}.  Let $M=B_{\bC}^2/\Sigma$ be
a $\bZ_2\times\bZ_2$ cover of $X$ as above.  From Lemma 1, we conclude that $h^{1,0}(M)=0$, from which we conclude from the discussions in
the proof of Lemma 4 that $h^0(M,K_M)=3$.   Hence the canonical map $\Phi$ is apriori a birational map from $M$ to $P_{\bC}^2$.
Lemma 1 also implies that the Picard number $\rho(M)=1$.  
From Proposition 1, we conclude that the canonical map is base point free and hence is a well-defined holomorphic map.
The degree of the canonical map is given by 
$$\int_M\Phi^*\cO(1)\cdot \Phi^*\cO(1)=K_M\cdot K_M=4K_X\cdot K_X=36,$$
since $X$ is a fake projective plane and hence $K_X\cdot K_X=9.$  
The surface is minimal since it is a complex ball quotient and hence does not contain
rational curves due to hyperbolicity of $M$.  Theorem 1 follows.  

\qed

\bigskip

\noindent{\bf References} 

\ms
\ni [B] Beauville, L'application canonique pour les surfaces de type g\'en\'eral, Inv. Math. 55(1979), 121-140.

\ms
\ni [BR] Blasius, D., Rogawski, J., Cohomology of congruence subgroups of  $SU(2,1)^p$ and Hodge cycles on some special complex hyperbolic surfaces. Regulators in analysis, geometry and number theory, 1-15, Birkh\"auser Boston, Boston, MA, 2000.

\ms
\ni [CV] Calabi, E., Vesentini, E., On compact locally symmetric K\"ahler manifolds, Ann. of Math. 71 (1960),
472-507.

\ms
\ni [CS] Cartwright, D., Steger, T.,  Enumeration of the $50$ fake projective planes, C. R. Acad. Sci. Paris, Ser. 1,
348 (2010), 11-13,  see also \\
\verb'http://www.maths.usyd.edu.au/u/donaldc/fakeprojectiveplanes/'

\medskip
\noindent
[DM]
Deligne, P., Mostow, G. D., Commensurabilities among lattices in $
PU(1,n)$. Annals of Mathematics Studies, 132. Princeton University Press, 
Princeton, NJ, 1993.

\ms
\ni [DPP] Deraux, M., Parker, J. R., Paupert, J., New non-arithmetic complex hyperbolic lattices, Invent. Math. 203 (2016), 681-771.

\ms
\ni
[DG] Du, R., and Gao, Y., Canonical maps of surfaces defined on abelian covers, Asian Jour. Math. 18(2014), 219-228.

\ms 
\ni [GKMS] Galkin, S.,  Katzarkov, L., Mellit A., Shinder, E., Derived categories of Keum's fake projective planes, Adv. Math. 278(2015), 238-253.

\ms
\ni [HL] Hambleton, I., Lee, Ronnie, Finite group actions on $P^2(\bC)$, Jour. Algebra 116(1988), 227-242.

\ms
\ni [H] Hirzebruch, F.,  Arrangements of lines and algebraic surfaces. Arithmetic and geometry, Vol. II, 113 -140, Progr. Math., 36, Birkh\"auser, Boston, Mass., 1983.

\ms
\ni [K] Keum, J, A vanishing theorem on fake projective planes with enough automorphisms, arXiv:1407.7632v1, Tran. Amer. Math. Soc. 369(2017), 7067-7083.

\ms
\ni [LY] Lai, C.-J., Yeung, S.-K., Exceptional collection of objects on some fake projective planes, IMRN, https://doi.org/10.1093/imrn/rnab186

\ms
\ni
[Mu] Mumford, D., An algebraic surface with $K$ ample, $K^2=9$, $p_g=q=0.$
 Amer. J.
Math. 101 (1979), 233--244.

\ms
\ni
[Pa] Pardini, R., Canonical images of surfaces, J. reine angew. Math.
417(1991), 215-219.

\ms
\ni [Pe] Persson, U., Double coverings and surfaces of general type., In: Olson,L.D.(ed.) Algebraic geometry.
(Lect. Notes Math., vol.732, pp.168-175) Berlin Heidelberg New York: Springer 1978.

\ms
\ni [Pr] Prasad, G., Volumes of S-arithmetic quotients of semi-simple groups. Publ. Math.,
Inst. Hautes \'Etud. Sci. 69, 91-117 (1989)

\ms
\ni [PY] Prasad, G., and Yeung, S.-K.,  Fake projective planes. Inv.\,Math. 168(2007), 321-370; Addendum, ibid 182(2010), 213-227.
%with a combined version in arXiv math/0512115.

\ms
\ni [R\'e] R\'emy, R., Covolume des groupes S-arithm\'etiques et faux plans projectifs, [d'apr\`es Mumford, Prasad, Klingler, Yeung,
Prasad-Yeung], S\'eminaire Bourbaki, 60\`eme ann\'ee, 2007-2008, no. 984.

\ms
\ni[Re] Reznikov, A., Simpson's theory and superrigidity of complex hyperbolic lattices, C. R. Acad.
Sci. Paris Sr. I Math., 320(1995), pp. 1061-1064.

\ms
\ni [Ro] Rogawski, J., Automorphic representations of the unitary group in three variables, Ann.
of Math. Studies, 123 (1990).

\ms
\ni [S] Su, J. C., Transformation groups on cohomology projective spaces, Trans. Amer. Math. Soc. 106 (1963), 305-318.

\ms
\ni [T] Tan, S.-L., Surfaces whose canonical maps are of odd degrees. Math. Ann.
292 (1992), 13-29.

\ms
\ni [W] Wilczy\'nski, D. M.,  Group actions on the complex projective plane, Trans. Amer. Math. Soc. 303 (1987), 707-731.

\ms
\ni [X] Xiao, G., Algebraic surfaces with high canonical degree, Math.
Ann., 274(1986), 473-483.

\ms
\ni [Ya] Yau, S.-T., Calabi's conjecture and some new results in algebraic geometry, Proc. Nat. Ac. Sc. USA 74(1977), 1798-1799.

\ms
\ni [Ye1] Yeung, S.-K., Integrality and arithmeticity of co-compact lattices corresponding to certain complex two
ball quotients of Picard number one, Asian J. Math. 8 (2004), 104-130; Erratum, Asian J. Math. 13
(2009), 283-286.

\ms
\ni [Ye2]  Yeung, S.-K., Classification and construction of fake projective planes, Handbook of geometric analysis, No. 2, 391-431, Adv. Lect. Math. (ALM), 13, Int. Press, Somerville, MA, 2010.

\newpage
\begin{center}
{\bf{\sc Corrigendum}}
\end{center} 

\bigskip
\noindent

\ni{\bf 1.} In the paper [Y], the proof of Lemma 6 contains the following error pointed out by Keum and F. Catanese (cf. arXiv:1801.05291).
In the middle of {\bf 7.2} for the proof of Lemma 6, there is the wrong claim that
 the bundle $K_X+\tau_i$ and hence $t_i$ is invariant under
$\bZ_7$ on $X$.  
The goal here is to give a completely new proof of Lemma 6 and hence Theorem 1.  All the unexplained notations are referred to [Y].

\ms
  Recall that  $X=B^2_{\bC}/\Pi$ is a fake projective plane with $\Pi$ one of the lattices
found in {\bf 5.11, A3} of [16], and is labelled as $(a=7,p=2,\emptyset,D_32_7)$ in [4]
as a surface in the class corresponding to $(a=7,p=2,\emptyset, D_32_7)$.  
The automorphism group of~$X$ is $G=\langle a,b\ |\ a^7=b^3=1, bab^{-1}=a^2\rangle$.
It contains just one subgroup $\langle a\rangle$ of order~7,
and the subgroup $\langle b\rangle$, whose seven distinct conjugates
are the Sylow~3 subgroups of~$G$.

%The automorphism group of $X$ is 
%$G=\bZ_7\rtimes \langle b\rangle=\langle a,b|a^7=b^3=1, bab^{-1}=a^2\rangle.$  $G$ contains a normal subgroup $\langle a\rangle$ of order $7$, and seven
%conjugate Sylow $3$ subgroups, one of which we denote by $G_3$.  
The surface $M$ that is to be proved to have canonical degree $36$ is constructed as follows.
It is known from
the file of registry of surface of [4] that $H_1(X/G,\bZ)=\bZ_2$, 
$H_1(X/\langle a\rangle,\bZ)=\bZ_2$, 
$H_1(X/\langle b\rangle,\bZ)=\bZ_2^2$ and $H_1(X,\bZ)=\bZ_2^4$.  Let $\rho_1:X\rightarrow X/\langle b\rangle$ and $\rho_2:X\rightarrow X/\langle a\rangle$ and $\rho:X\rightarrow X/G$ 
be the
projection maps.
Denote by $A\cong \bZ_2^2$ the first two factors of $\bZ_2$ in $H_1(X,\bZ)$ so that $(\rho_1)_*(A)=H_1(X/G_3,\bZ)$, and $B$ the first factor of $\bZ_2$  in $A$ 
so that $\rho_*(A)=H_1(X/G,\bZ)$.  In such case, $(\rho_2)_*(B)=H_1(X/\langle a\rangle,\bZ)$.  Let $p:M\rightarrow X$ be the $\bZ_2\times\bZ_2$ cover of $X$ with
fundamental group $\Sigma$ obtained by
kernel of the homomorphism of $\alpha:\Pi\rightarrow A$.  In other words, $\Pi/\Sigma =A.$   From construction,
$p^*(K_X+\tau_i)=K_M$ for $i=1,2,3$, since $\tau_i$'s corresponds to elements of $A$, and $M$ is a covering of $X$ given by $\ker\alpha$.

%We remark that the same construction works for any of the Sylow $3$-subgroup of $G$.  However as the Sylow  $3$ subgroups 
%are all
%conjugate to one another by an element in $\langle a\rangle$, there would only be one $M$ constructed up to conjugation by an element in 
%$\Lambda$ from the seven Sylow 3 subgroups.

\ms
The argument before
{\bf 7.2} of [Y] is valid.  In particular, the following is proved, \\
(i) from {\bf 4.3} of [Y], $\Gamma(M,K_M)$ has dimension $3$ and from the second paragraph of {\bf 7.1}, is spanned by $s_1, s_2, s_3 \in \Gamma(M,K_M)$, where the
zero divisor $Z_{s_i}=\{s_i=0\}$ of $s_i$ is invariant as a set under $A$;\\
(ii) from the first paragraph of {\bf 7.2}, $s_i$ descends to a section $t_i\in \Gamma(X,K_X+\tau_i)$ for some two torsion line bundle $\tau_i$, which could be identified with a non-zero element in $A\cong \bZ_2^2$ 
by the Universal Coefficient Theorem for $i=1,2,3$, in the sense that the torsion line bundles can be identified with the torsion elements in $H_1( X,\bZ)$;\\
(iii) $h^0(X,K_X+\tau_i)=1$ for all $i=1,2,3$ since $h^0(M,K_M)=3$ as mentioned above, and as $t_i\in H^0(X,K_X+\tau_i)$, we conclude that
 $t_i^2\in \Gamma(X,2K_X)$;\\ 
(iv) from Lemma 4 and 5 of [Y], there is no base locus of $\Gamma(M,K_M)$ except possibly at a finite number of points, that is, $\cap_{i=1}^3Z_{s_i}$ contains at most a finite number of points,
where $Z_s$ denotes the zero divisor of $s$.

\ms
For Theorem 1, it suffices
to prove Lemma 6 in \S7 of [Y] in the sense that such isolated based points do not exist.   In this Corrigendum, this is shown by relating sections of $\Gamma(M,K_M)$ to certain sections
of $\Gamma(X,2K_X)$.
%In other words, $\cap_{i=1}^3Z_{s_i}$ does not contain isolated points for sections $s_i\in \Gamma(M,K_M)$ chosen as in [Y],
%each of which invariant as a set under the covering group $\bZ_2\times \bZ_2$ of $M\rightarrow X$ and hence descends to $X$.  
%In this way, $M$ provides an example to
%Theorem 1.  
We remark that $M$ as described above
has a presentation given by {\bf 3.3} of [Y] after checking with Magma, though the explicit presentation is not needed for the proof of Theorem 1 here.

\ms
\ni{\bf 2.}  The automorphism group~$G$ has three 1-dimensional representations, $\chi_0$,
$\chi$ and~$\bar\chi$, and two 3-dimensional irreducible representations $\pi$, $\bar\pi$.
Here $\chi_0$ is the trivial character, and $\chi(a)=1$ and $\chi(b)=\zeta_3$, 
while $\pi(a)$ is the diagonal matrix with diagonal entries 
$\zeta_7$, $\zeta_7^2$ and~$\zeta_7^4$, $\pi(b)$ is the permutation matrix 
corresponding to the permutation $(1,3,2)$, and $\bar\chi$ and $\bar\pi$ 
are the complex conjugates of~$\chi$ and~$\pi$, and $\zeta_n$ is a fixed
primitive $n$-th root of unity.

%It is a standard fact that the automorphism group $G$ has character table given by

%$$\hskip-0.8in
%\begin{array}{c|ccccc}
%&1&a&a^3&b&b^2\\
%\hline
%\chi_1&1&1&1&1&1\\
%\chi_2&1&1&1&\zeta_3&\zeta_3^2\\
%\chi_3&1&1&1&\zeta_3^2&\zeta_3\\
%\chi_4&3&\alpha&\beta&0&0\\
%\chi_5&3&\beta&\alpha&0&0
%\end{array}
%$$
%where $\alpha=\zeta_7+\zeta_7^2+\zeta_7^4, \beta=\zeta_7^3+\zeta_7^5+\zeta_7^6$ and $\zeta_n$ is a fixed primitive $n$-th root of unity.

Consider the $G$-spaces $V=H^0(X,2K_X)$ and $P_{\bC}^9=P_{\bC}(V)$. They contain
one copy of the trivial representation, 2~copies of~$\pi$ and 1~copy of~$\bar\pi$, so that $V=V_0+2V_1+V_2$.
This is given explicitly in [BK] (2.1): % Perhaps you are using~\cite[(2.1)]{BK}:
\begin{eqnarray}
&&a(u_0:\ u_1:u_2:u_3:u_4:u_5:u_6:u_7:u_8:u_9)\notag\\
&=&(u_0:\zeta_7^6u_1:\zeta_7^5u_2:\zeta_7^3u_3:\zeta_7u_4:\zeta_7^2u_5:\zeta_7^4u_6:\zeta_7u_7:\zeta_7^2u_8:\zeta_7^4u_9)\\
&&b(u_0:\ u_1:u_2:u_3:u_4:u_5:u_6:u_7:u_8:u_9)\notag\\
&=&(u_0:u_2:u_3:u_1:u_5:u_6:u_4:u_8:u_9:u_7)
\end{eqnarray}
%Hence we may decompose into irreducible representations, where $V_1$ is dimension one spanned by $u_0$, which is
%a section of $H^0(X,2K_X)$ invariant under $G$ corresponding to $\chi_0$, and $V_1$ (resp. $V_2$) are irreducible representations corresponding to %$\rho$ (resp. $\bar\rho$), where $\rho$ is either 
%$\pi$ (resp. $\bar\pi$).
 
We know from Riemann-Roch and Kodaira Vanishing Theorem that $h^0(X,2K_X)=10$.   Denote by $h^0(X,2K_X)^{H}=\dim_{\bC}\Gamma(X,2K_X)^H$
the dimension of the subspace of sections
invariant up to a scalar multiple under a group $H$.
We find that $h^0(X,2K_X)^{\langle b\rangle}=4$.  To see this,
 from [4] or [11], we know that the singular set of $X/\langle b\rangle$ consists of three $\frac13(1,2)$ points,
the resolution of each is a chain of two $(-2)$ curves.  Hence if $\sigma:Y\rightarrow X/\langle b\rangle$ is the canonical resolution, $K_Y=\sigma^*K_{X/\langle b\rangle}$ so that
$K_Y^2=(\sigma^*K_{X/\langle b\rangle})^2=3$ and $c_2(Y)=9$, which gives rise to $h^0(X,2K_X)^{\langle b\rangle}=h^0(X/\langle b\rangle,2K_{X/\langle b\rangle})=4$ from Riemann-Roch and Kawamata-Viehweg
Vanishing Theorem.  Alternatively, this also follows from the representation above.  Furthermore, 
the representation above shows that 
$h^0(X,2K_X)^G=h^0(X/G,2K_{X/G})=1$, where $K_{X/G}$ is regarded as a $\bQ$ line bundle.  
%Note that from [4] and [11], $X/G$ has
%three singularities of type $\frac13(1,2)$ which is resolved to a chain of two $(-2)$ curves, and a singularity of type $\frac17(1,3)$ which resolves
%to a chain $C-E-F$ where $C$ is a $(-3$) curve, and $E, F$ are $-2$ curves.  Let $\pi:Y\rightarrow X/G$ be the canonical resolution of $X/G$.  It follows from intersection with
%various exceptional curves that  
%$K_W=\pi^*K_Y-\frac37 C-\frac 27 E-\frac17 F.$ It follows from direct computations that $c_1^2=0$ and $c_2=12$.  Hence $\chi(\cO)=1$ from Noether's Formula.

We would use (i)-(iv) in {\bf 1}, the explicit computations of [BK], and study of $\langle b\rangle$ invariant sections to check that there is no isolated points in $\cap_{i=1}^3Z_{s_i}.$  We used Magma which is symbolic and exact.

\ms
\ni{\bf 3.}  {\bf Proof of Lemma 6}  As summarized in (i)-(iv) above, $H^0(M,K_M)$ is generated by $s_i, i=1,2,3$, which descend to effective sections $t_i$ of $K_X+\tau_i$ on $X$. 
% And as $h^0(M,K)=3$, $h^0(X,K+\tau_i)=1$ for all $i=1,2,3$.  
From our setting in {\bf 1}, each $\tau_i, i=1,2,3$ is invariant under $\langle b\rangle$ corresponding to the non-trivial element
in $H_1(X/\langle b\rangle,\bZ)=\bZ_2^2$ from the Universal Coefficient Theorem.  One of those three, say denoted by $\tau_3$, is invariant
under $G$ corresponding to $H_1(X/G,\bZ)=\bZ_2$.  Hence under the action of $\langle a\rangle$, the torsion line bundle $\tau_3$ is fixed,
while $\tau_1, \tau_2$ are not invariant under $\langle a\rangle$ and hence $a$ acts freely within each orbit $\langle a\rangle\tau_1$ and 
$\langle a\rangle\tau_2$. 
The two orbits are disjoint, for if $a^k\tau_1=\tau_2$ for some $1\leqslant k\leqslant 6$, then
$$\tau_2=b\tau_2=ba^k\tau_1=ba^kb^{-1}\tau_1=a^{2k}\tau_1=a^k\tau_2,$$
contradicting the free action of $\langle a\rangle$ on the orbit of $\tau_2$.
It follows that the cardinality of the set $\cup_{i=1}^2\langle a\rangle \tau_i$ is $14$, which together with the trivial bundle $0$ and $\tau_3$ exhaust 
the torsion line bundles of $X$ corresponding to $H_1(X,\bZ)=\bZ_2^4$ from the Universal Coefficient Theorem.  In addition to $\tau_3$, we denote the remaining
fourteen non-trivial $2$-torsion line
bundles by $\sigma_j, j=1,\dots,14$.

Since $h^0(X,K_X+\tau_i)=1$ for $i=1,2,3$,  it follows under the action of $\langle a\rangle$ as explained above that 
$h^0(X,K_X+\tau)=1$ for all $2$-torsion line bundle $\tau\neq 0$ corresponding to $H_1(X,\bZ)=\bZ_2^4$.  
Hence 
there are $14$ sections $w_j, 1\leq j\leq 14$ in $\cup_{i=1}^{14}\Gamma(X,K_X+\sigma_i)$ corresponding to
$H_1(X,\bZ)=\bZ_2^4.$  The $14$ divisors consist of two $\langle a\rangle$ orbits, $\langle a\rangle t_2$ and $ \langle a\rangle t_3$.   
 The square of each such section gives rise to a section of $V=\Gamma(X,2K_X)$.  
%We claim that $t_1^2\in W_1$ is not contained in $V_2$.  
From (2), the vector space 
$\Gamma(X,2K_X)^{\langle b\rangle}\cong\Gamma(X/\langle b\rangle, 2K_{X/\langle b\rangle})$ has a basis given explicitly from (2) by $v_0:=u_0$,
 $v_1:=u_1+u_2+u_3$, $v_2:=u_4+u_5+u_6$, $v_3:=u_7+u_8+u_9$.

 We already know that $t_1^2=u_0$ from [BK], which also follows from the command {\verb IsDomain  in Magma.   It suffices for us to show that $\cap_{j=1}^3Z_{t_i}=\emptyset$, for which we would give two arguments.
 
 The first proof is to use reduction at a finite field $F_p$,  where $p$ is chosen to be $23$ for convenience, and utilizing comparsion theorem of Grothendieck as given
 in SGA, XII 7, [G].
 By checking the three $7\times 7$ minors of Jacobians of the defining functions of $X$ given in [BK] using Magma, we verify that $X$ is smooth for $p=23$ and hence has good reduction
 at $23$, which we denote by $X^p$.   As $\pi_1(X)$ is residually finite, we can identify the topological fundamental group with its etale fundamental group.
 The first homology group $H_1(X,\bZ)$, as $\pi_1(X)$ modulo its commutator, is identified with abelianization of the maximal quotient of the etale fundamental group by a prime relatively prime to $p$.
 The same principle holds for the resolution of $X/\langle b\rangle$ at its three singular points and hence for $X/\langle b\rangle$.
 As $p\neq 2$ and $H_1(X,\bZ)^{\langle b\rangle}=\bZ_2^2$, we conclude that $H_1(X^p,\bZ)^{\langle b\rangle}=\bZ_2^2$.  
% From Correspondence Principle, we check that there are precisely three non-zero $2$-torsion line bundles $\tau_i, i=1.2.3$
 % in $M^p$ invariant under
% $\langle a\rangle$ corresponding to $H^1(X^p,\bZ)^{(\langle b\rangle)}=\bZ_2^2$, and $h^0(X^p,2K_{X_p})^{(\langle b\rangle)}=3$, with three linearly independent sections
 % corresponding to square of a section 
% $s_{p,i}\in H^0(X^p,K_{X^p}+\tau_i)$,  each of which has dimension $1$.  Here e'tale cohomology is used and $p$ is relative prime to $2$.  
 In particular, a non-trivial $2$-torsion line bundle on $M$ gives rise to a non-trivial $2$-torsion line
 bundle on $M_p$.  This implies that $t_i, i=1,2,3$ would give $3$ different image $t_{p,i}$ on reduction modulo $p$.  
 %In general under reduction modulo a prime $p$, 
% we can only conclude that $h^0(X^p,2K_{X^p})^{\langle b\rangle}\geqslant h^0(X,2K_{X})^{\langle b\rangle}$ and $h^0(X^p,K_{X^p}+\tau)^{\langle b\rangle}\geqslant h^0(X,K_{X}+\tau)^{\langle b\rangle}$ for a $2$-torsion line bundle $\tau$.
 Recall that sections of $H^0(X,2K_{X})^{\langle b\rangle}$ are spanned by square of sections of $H^0(X,K_X+\epsilon)^{\langle b\rangle}$ for some $2$-torsion line bundle $\epsilon$.
From earlier discussion, $H^0(X^p,2K_{X^p})^{\langle b\rangle}$ has dimension $4$ and $t_i^2, i=1,2,3$ are linear combinations of $v_j, j=1,\dots,4$.
It follows that $t^2_{p,i}$ has to be a linear combination of $v_j$ and is reducible or non-reduced modulo $p$. 
 We 
 apply \verb IsDomain  in Magma to each section $\sum_{i=0}^3c_iv_i$ for $c_i\in \{0,\dots,22\}$ and find that there are exactly three quadruples $c_i$ for which the sections are reducible or non-reduced, given by $v_0, v_0+14v_1$ and  $v_0+22v_1+11v_2+19v_3$.
Since we know already that there are three such sections coming from $t^2_{p,i}, i=1,2,3$, these have to be $t^2_{p,i}$.   We check that $\cap_{i=1}^3 Z_{t^2_{p,i }}=\emptyset$ on $X_p$ from the command \verb HilbertPolynomial  in Magma, which gives value $0$.  This implies that $\cap_{j=1}^3Z_{t^2_i}=\emptyset$ on $X$,
which leads to  $\cap_{j=1}^3Z_{t_i}=\emptyset$ 
on $X$ and hence $\cap_{j=1}^3Z_{s_i}=\emptyset$ 
on $M$.

 The second proof is more explicit.  Recall that $X$ is a Shimura variety and is defined over a number field $\bQ(\sqrt{-7})$.  Recall that $t^2_1=u_0$.
 Since $t_2, t_3$ cannot be deformed as curves
 from the fact that $h^0(X,K_X+\tau_i)=1$ for all $2$-torsion $\tau_i$, we know that they are rigid and can be defined over $\overline\bQ$.  Let $j=2, 3$.  Since the ring of
 integers $\cO_{\bQ(\sqrt {-7})}$ has a basis given by $1,\eta=\frac12(1+\sqrt{-7})$,
 we try $t^2_j=v_0+\sum_{j=2}^3 (\alpha_j+\beta_j \eta)v_j$ for some $\alpha_j,
 \beta_j\in\bQ$. 
 By considering  reduction  modulo $p=11, 23, 29$ and using \verb IsDomain  as in the first method, we conclude that a candidate for $t_2$ is $\hht_2$ with $\hht_2^2=v_0+\eta v_1$,
 corresponding to $v_0+14v_1$ for $p=23$.
 Using  \verb IsDomain  command over $\bQ(\sqrt{-7})$, we conclude that $\hht_2$ is either non-reduced or reducible and hence
 has to be square of some section of a bundle numerically equivalent to $K_X$ from proof of Lemma 2 of [Y].  Hence  $t_2=\hht_2$ up to a scaling constant.
 Similar procedure leads to $t_3^2=v_0+(-6+2\eta)v_1+(8-8\eta)v_2-4v_3$  up to a scalar.
 %We can now complete the argument in one of the
 %following arguments.  The first argument is to use a similar argument to show that $s_3=u_0+(-6+2\eta)v_1+8(1-\eta)v_2-4v_3$ is another candidate.   Using \verb HilbertPolynomial, we check that $\cap_{i=1}^3\bZ_{s_i }=\emptyset$.
 %The second argument to use the fact that 
Since  the image of $t_3$ in $X^p$ in reduction modulo $23$ is $t_{p,3}$ studied earlier, and $t_{p,3}$
does not have non-trivial intersection with the intersection of $Z_{t^2_1}\cap Z_{t^2_2}$ modulo $p=23$ from command \verb HilbertPolynomial  in Magma, 
we conclude that $\cap_{j=1}^3Z_{t^2_i}=\emptyset$.  Alternatively, we show that $Z_{t^2_1}\cap Z_{t^2_2}$ actually occurs only at explicit points given
 by the three fixed points of $\bZ_3$ on $X$.  Using \verb HilbertPolynomial, one shows that $Z_{v_3}$ does not intersect $Z_{t_1}\cap Z_{t_2}=Z_{v_0}\cap Z_{v_1}$ but $Z_{v_2}$ does.  Hence
 if $t^2_3=\sum_{i=0}^3c_iv_i$ and
  $\cap_{i=1}^3 Z_{t_i }\neq\emptyset$, the only possibility is that $c_3=0$ by evaluating $t_3^2$ at $\cap_{j=0}^2Z_{v_j}$.
 This contradicts the earlier fact that $t^2_{p,3}= v_0+22v_1+11v_2+19v_3$, corresponding to $c_3=19 \pmod{23}$.

% We remark that Carlos Rito has mentioned to the author that he had also a proof of the emptyness of $\cap_{i=1}^3Z_{s_i}$.
 % Note that $-\eta=\zeta_7^6+\zeta_7^5+\zeta_7^3$ and $-\bar\eta=\zeta_7+\zeta_7^2+\zeta_7^4$.
 
 \bs
\noindent{\bf References} 

\ms
\ni [BK] Borisov, L. A., Keum, J., Explicit equations of a fake projective plane,\\
 arXiv:1802.06333, Duke Math. J. 169(2020), 1135-1162. 

\ms
\ni [G] Grothendieck, A., SGA1, arXiv:math/0206203v2.

\ms
\ni [Y] Yeung, S.-K., A surface of maximal canonical degree,  Math. Ann. 368(2017), 1171-1189.

\end{document}